\documentclass[10pt]{article}

\RequirePackage{amsmath}
\RequirePackage{amsfonts}
\RequirePackage{amssymb}
\RequirePackage{hyperref}
\RequirePackage{multirow}
\RequirePackage{algorithm}
\RequirePackage{algpseudocode}
\RequirePackage{graphicx}
\RequirePackage{xcolor}
\RequirePackage{amsthm}
\RequirePackage{amsmath}
\RequirePackage{amsfonts}
\RequirePackage{amssymb}
\RequirePackage{hyperref}
\RequirePackage{multirow}
\RequirePackage{algorithm}
\RequirePackage{algpseudocode}
\RequirePackage{graphicx}
\RequirePackage{xcolor}
\allowdisplaybreaks
\usepackage[sort&compress]{natbib}

\usepackage{dpr_fec,ifthen,dsfont}

\newtheorem{assumption}{Assumption}[section]

\newtheorem{theorem}{Theorem}[section]

\usepackage{isetechreport}
\def\reportyear{26T}
\def\reportno{001}
\def\revisionno{0}
\def\originaldate{Jan 16, 2026}
\def\revisiondate{\today}
\coraltrue
\cvcrfalse

\begin{document}

\title{Projected Stochastic Momentum Methods for Nonlinear Equality-Constrained Optimization for Machine Learning}

\author{Qi Wang}
\affil{Department of Industrial and Operations Engineering, University of Michigan\thanks{qiwangqi@umich.edu}}
\author{Christian Piermarini}
\affil{Amazon S\`{a}rl\thanks{chripiermarini@gmail.com}}
\author{Yunlang Zhu, Frank E.~Curtis}
\affil{Department of Industrial and Systems Engineering, Lehigh University\thanks{yuza23@lehigh.edu}\thanks{frank.e.curtis@lehigh.edu}}

\titlepage

\maketitle

\begin{abstract}
  Two algorithms are proposed, analyzed, and tested for solving continuous optimization problems with nonlinear equality constraints.  Each is an extension of a stochastic momentum-based method from the unconstrained setting to the setting of a stochastic Newton-SQP-type algorithm for solving equality-constrained problems.  One is an extension of the heavy-ball method and the other is an extension of the Adam optimization method.  Convergence guarantees for the algorithms for the constrained setting are provided that are on par with state-of-the-art guarantees for their unconstrained counterparts.  A critical feature of each extension is that the momentum terms are implemented with projected gradient estimates, rather than with the gradient estimates themselves.  The significant practical effect of this choice is seen in an extensive set of numerical experiments on solving informed supervised machine learning problems.  These experiments also show benefits of employing a constrained approach to supervised machine learning rather than a typical regularization-based approach.
\end{abstract}

\newcommand{\true}{\text{true}}
\newcommand{\edit}[1]{\textcolor{blue}{#1}}

\section{Introduction}

Algorithms based on the stochastic-gradient methodology \cite{RobbMonr51,RobbSieg71} have been found to be especially powerful for solving modern-day unconstrained continuous optimization problems that arise in multiple areas, most notably in supervised machine learning.  Of the many variants of the stochastic-gradient methodology, momentum-based approaches have been particularly popular and effective in practice.  These include the heavy-ball method~\citep{polyak1964some}, Adagrad~\citep{DuchHazaSing11}, RMSprop~\citep{DaupVrie15,TielHint12}, and Adam \citep{KingBa14}; see \citep{BottCurtNoce2018} for an overview.

The main contributions of this paper are extensions of the heavy-ball and Adam approaches from the unconstrained setting to the setting of a stochastic Newton-SQP framework for solving nonlinear-equality-constrained continuous optimization problems.  Our approaches are inspired by the stochastic Newton-SQP method proposed in~\citep{BeraCurtRobiZhou21} (see also \citep{BeraCurtONeiRobi23,CurtONeiRobi23}) for solving problems with nonlinear equality constraints.  We show that our proposed methods can offer theoretical guarantees that are on par with state-of-the-art guarantees that have been offered for the heavy-ball method and Adam in the unconstrained setting.  In particular, our analyses follow the analyses for these methods in the unconstrained setting that are presented in \citep{defossezSimpleConvergenceProof2022}.

Due to its more impressive practical performance, our more significant contribution in this paper is the extension of Adam to the equality-constrained setting.  A different extension of Adam to the equality-constrained setting (without a convergence guarantee) has been proposed previously; see \citep{Marquez-Neila}.  However, besides the fact that we offer a theoretical convergence guarantee for it, our algorithm is unique in that the running averages that are maintained in an Adam-based approach are taken with \emph{projected} gradient estimates rather than with gradient estimates themselves.  This is also the case with our extension of the heavy-ball method.  We show in an extensive set of numerical experiments with informed supervised machine learning test problems that our projected stochastic Adam algorithm can outperform both (a) Adam applied to minimize a regularized objective function and (b) the projection-less extension of Adam proposed in \citep{Marquez-Neila}.

A variant of the stochastic Newton-SQP method from \cite{BeraCurtRobiZhou21} that employs Adagrad-type scaling has been proposed and analyzed in \cite{ONeil2024}.  However, the algorithm and analysis in this paper are distinct from those in~\cite{ONeil2024} due to the subtle, yet significant differences between convergence analyses for Adagrad- versus other methods.  As far as we are aware, ours is the first paper that offers convergence guarantees for heavy-ball- and Adam-based stochastic Newton-SQP methods for solving nonlinear-equality-constrained problems.

\subsection{Outline}

The class of optimization problems of our interest and fundamental properties of the algorithms that we propose and analyze are presented in~\S\ref{sec.problem_algorithm}.  Our proposed heavy-ball- and Adam-based schemes, along with our convergence analyses of them, are given in \S\ref{sec.analysis}.  In \S\ref{sec.implementation}, we discuss a broad class of problems, namely, informed supervised machine learning problems, for which the proposed algorithms are particularly well suited.  In that section we also discuss a few critical considerations for efficient implementations of the algorithms.  In \S\ref{sec.experiments}, we present the results of a large set of experiments that demonstrate the effectiveness of the algorithms.  Finally, we provide some concluding remarks in \S\ref{sec.conclusion}.

\section{Problem Formulation and Stochastic Newton-SQP Framework}\label{sec.problem_algorithm}

Our problem class of interest is that of continuous equality-constrained optimization problems, where the objective function $f : \R{n} \to \R{}$ and constraint function $c : \R{n} \to \R{m}$ are continuously differentiable.  Our proposed algorithms are designed to solve such problems when the objective function is defined by an expectation of a function with a random variable argument and may be nonconvex, and when the equality-constraint function may be nonlinear.  Formally, our problem class of interest can be expressed as instances of
\bequation\label{prob.opt}
  \min_{x \in \R{n}}\ f(x)\ \ \text{subject to}\ \ c(x) = 0,\ \ \text{where}\ \ f(x) = \E_\xi[F(x, \xi)]\ \ \text{for all}\ \ x \in \R{n},
\eequation
$\xi$ is a random variable with associated probability space $(\Xi,\Fcal_\xi,\P_\xi)$, $F : \R{n} \times \Xi \to \R{}$, and $\E_\xi$ denotes the expected value operator with respect to the probability measure $\P_\xi$.  At a given $x \in \R{n}$, a first-order optimality condition for problem \eqref{prob.opt} is that there exists a Lagrange multiplier $y \in \R{m}$ such that
\bequation\label{eq.firstorder_y}
  \nabla f(x) + J(x)^T y = 0\ \ \text{and}\ \ c(x) = 0,\ \ \text{where}\ \ J := \nabla c^T.
\eequation
The first of these conditions states that the gradient $\nabla f(x)$ lies in the range space of the constraint derivative~$J(x)^T$.  Due to the fundamental theorem of linear algebra, this is equivalent to the property that the projection of the gradient $\nabla f(x)$ onto the null space of the constraint Jacobian $J(x)$ is equal to the zero vector.  Under the assumption that the constraint Jacobian has full row rank and with $P(x)$ denoting the projection operator onto this null space at $x$, the first-order optimality conditions in \eqref{eq.firstorder_y} are equivalent to
\bequation\label{eq.firstorder}
  P(x) \nabla f(x) = 0\ \ \text{and}\ \ c(x) = 0,\ \ \text{where}\ \ P(x) = I - J(x)^T (J(x) J(x)^T)^{-1} J(x).
\eequation
These are the form of the first-order optimality conditions that we employ in our analyses.

The algorithms that we propose, analyze, and test in this paper have as a basis the stochastic-gradient-based Sequential Quadratic Programming (SQP) framework proposed and analyzed in \citep{BeraCurtRobiZhou21} for solving equality-constrained optimization problems.  A simplified version of this method is stated as Algorithm~\ref{alg.sqp}.  The key aspect of it that is distinct for the setting of constrained optimization is that each search direction is computed to satisfy a linearization of the constraints that is defined with respect to the current iterate $x_k$.  Specifically, given an objective gradient estimate $g_k \in \R{n}$ and a symmetric and positive-definite matrix $H_k \in \R{n \times n}$, the search direction $d_k$ can equivalently be defined as the solution of the quadratic optimization subproblem
\bequation\label{eq.qp}
  \min_{d \in \R{n}}\ g_k^Td + \thalf d^TH_kd\ \st\ c(x_k) + J(x_k)d = 0.
\eequation
The linear system \eqref{eq.linear_system} represents the necessary and sufficient conditions for optimality for \eqref{eq.qp}, so any solution of \eqref{eq.linear_system} yields a solution of \eqref{eq.qp} as well as a Lagrange multiplier estimate $y_k$.  (If $x \equiv x_k$ is a point at which there exists $y$ such that \eqref{eq.firstorder_y} holds and $g_k = \nabla f(x_k)$, then a solution of \eqref{eq.linear_system} is $(d_k,y_k) = (0,y)$.) It is well known that if $J(x_k)$ has full row rank and since $H_k$ is positive definite in the null space of $J(x_k)$, subproblem~\eqref{eq.qp} is feasible and has a unique globally optimal solution, which is given by the unique solution of~\eqref{eq.linear_system}.

\begin{algorithm}[ht]
  \caption{Stochastic-Gradient-based SQP Framework \citep{BeraCurtRobiZhou21}}
  \label{alg.sqp}
  \begin{algorithmic}[1]
    \Require $x_1 \in \R{n}$ and $\{\alpha_k\} \subset (0,1]$
    \For{all $k\in\N{}$}
      \State compute a stochastic gradient estimate $g_k \approx \nabla f(x_k)$ and choose symmetric $H_k \in \R{n \times n}$
      \State compute $d_k$ by solving
      \bequation\label{eq.linear_system}
        \bbmatrix H_k & J(x_k)^T \\ J(x_k) & 0 \ebmatrix \bbmatrix d_k \\ y_k \ebmatrix = - \bbmatrix g_k \\ c(x_k) \ebmatrix
      \eequation
      \State set $x_{k+1} \gets x_k + \alpha_k d_k$
    \EndFor
  \end{algorithmic}
\end{algorithm}

An algorithm that computes each search direction by solving a system of the form \eqref{eq.linear_system} can employ a direct solver for symmetric indefinite linear systems.  Alternatively, the search direction can be computed using a step decomposition method.  This is the approach that we specify for the methods that we propose in this paper since the methods make use of running averages of projected gradient estimates.  We close this section with a description of the linear algebra involved in a step decomposition approach, which also reveals the computational cost required to compute search directions in our proposed methods.  Importantly, these costs will be small when the number of equality constraints $m$ is small, meaning that in such cases the computational cost of each iteration will be proportional to the costs in an unconstrained context.

Suppose that, similarly as in Algorithm~\ref{alg.sqp}, an algorithm computes each search direction through solving
\bequation\label{eq.linsys}
  \bbmatrix I & J^T \\ J & 0 \ebmatrix \bbmatrix s \\ y \ebmatrix = - \bbmatrix q \\ c \ebmatrix,
\eequation
where $I$ is the identity matrix and we assume that $J \in \R{m \times n}$ has full row rank.  By the fundamental theorem of linear algebra, the solution component $s$ can be expressed as $s = v + u$, where $v \in \Range(J^T)$ and $u \in \Null(J)$.  Considering the second block of equations, one finds that $v$ can be computed by solving an $m$-dimensional positive-definite system $JJ^T\vtilde = -c$ for $\vtilde \in \R{m}$, then computing $J^T\vtilde = v$, which yields
\bequation\label{eq.v}
  v = - J^T(JJ^T)^{-1} c.
\eequation
The computational cost of computing $v$ is thus $\Ocal(m^3 + m^2 n)$.

Now letting $Z \in \R{n \times (n-m)}$ denote an orthogonal matrix whose columns span $\Null(J)$, the first row of \eqref{eq.linsys} states $u + J^Ty = -(q + v)$, so $u = -Z(Z^TZ)^{-1}Z^Tq$.  However, it is not efficient to compute~$u$ in this manner since it requires computing the null-space basis matrix $Z$.  Fortunately, one can replace $Z(Z^TZ)^{-1}Z^T$ with a matrix expressed in terms of $J$.  Specifically, one finds that $Z(Z^TZ)^{-1}Z^T = I - J^T (JJ^T)^{-1} J =: P$, so
\bequation\label{eq.u}
  u = -(I - J^T (JJ^T)^{-1} J) q = -Pq.
\eequation
In other words, $u$ is the negative of the projection of $q$ onto the null space of $J$; recall \eqref{eq.firstorder}.  Practically speaking, the component $u$ can be computed by computing the matrix-vector product $\qhat := Jq$, solving the $m$-dimensional positive definite system $JJ^T \qbar = \qhat$ for $\qbar \in \R{m}$, computing the matrix-vector product $J^T\qbar$, and adding the result to $-q$.  Thus, similar to $v$, the cost of computing $u$ is $\Ocal(m^3 + m^2 n)$.

Through this discussion, one can observe that the solution of \eqref{eq.linsys} is the same if $q$ is replaced by $Pq$.  This follows since \eqref{eq.v} reveals that $q$ has no effect on $v$, and since by virtue of $P$ being an orthogonal projection matrix one has by \eqref{eq.u} that $u = -Pq = -P^2q$.  That being said, the solution of the system is clearly affected if~$q$ is replaced by another vector not necessarily in the null space of the constraint Jacobian.  These comments reveal a critical distinction between our proposed version of Adam and that proposed in \citep{Marquez-Neila}.  In \citep{Marquez-Neila}, running averages (i.e., the momentum terms) are taken with gradient vectors, whereas in our proposed approach these running averages are taken with projected gradients.  This difference is critical theoretically and practically.

\section{Stochastic Momentum-based Algorithms}\label{sec.analysis}

In this section, we propose and analyze two new stochastic momentum-based methods for solving \eqref{prob.opt}.  Each algorithm computes search direction components through the formulas~\eqref{eq.v} and \eqref{eq.u}.  However, they are each distinct in the manner that scaling and momentum are applied to the latter component when constructing the search direction taken by the algorithm.

The stochastic nature of the algorithms means that our analysis of each of them considers a stochastic process defined by each algorithm.  In each case, let $(\Omega,\Fcal,\P)$ denote a probability space that captures the behavior of the algorithm, which is to say that each outcome in $\Omega$ represents a possible realization of a run of the algorithm.  In addition, let $\E$ denote the expected value operator defined by the probability measure $\P$.  The only source of randomness in each iteration is the computation of a stochastic gradient estimate.  Hence, using the notation of \eqref{prob.opt}, one can consider the set of outcomes as $\Omega \equiv \Xi \times \Xi \times \cdots$.  Let $\Fcal_1$ be the $\sigma$-algebra defined by the initial conditions of an algorithm, and, more generally, for all $k \in \N{}$ let $\Fcal_k$ be the $\sigma$-algebra generated by the initial conditions and the stochastic gradient estimators up through the end of iteration $k-1$.  In this manner, one has that $\Fcal_1 \subseteq \Fcal_2 \subseteq \cdots \subseteq \Fcal$ and the sequence $\{\Fcal_k\}$ is a filtration.

Throughout this section, we employ the shorthand notation $c_k := c(x_k)$, $J_k := \nabla c(x_k)^T$, and $P_k := P(x_k)$.

\subsection{Projected Stochastic Heavy-Ball SQP}

We first consider an extension of the heavy-ball method to the setting of equality-constrained optimization.  Specifically, we propose Algorithm~\ref{alg.heavy-ball}.

\begin{algorithm}[ht]
  \caption{Projected Stochastic Heavy-ball SQP}
  \label{alg.heavy-ball}
  \begin{algorithmic}[1]
    \Require $x_1 \in \R{n}$, $\{\rho_k\}$ and $\{h_k\}$ with $\rho_k \in (0,1]$ and $h_k \in \R{}_{>0}$ for all $k \in \N{}$, $\beta \in [0,1)$, and $\alpha \in (0,1]$
    \State set $r_0 \gets 0 \in \R{n}$
    \For{all $k\in\N{}$}
      \State compute a stochastic gradient estimate $g_k \approx \nabla f(x_k)$
      \State compute $v_k \gets - \rho_k J_k^T \left(J_k J_k^T\right)^{-1} c_k$
      \State compute $u_k \gets - h_k^{-1} P_k g_k$, where $P_k:= I - J_k^T (J_k J_k^T)^{-1} J_k$
      \State set $r_k \gets \beta r_{k-1} + u_k$
      \State set $d_k \gets v_k + P_k r_k$
      \State set $x_{k+1} \gets x_k + \alpha d_k$
    \EndFor
  \end{algorithmic}
\end{algorithm}

For our analysis of Algorithm~\ref{alg.heavy-ball}, we make the following assumption.

\begin{assumption}\label{ass.heavy-ball}
  There exists an open convex set $\Xcal \subseteq \R{n}$ containing the iterates of any run of Algorithm~\ref{alg.heavy-ball} over which the objective function $f : \R{n} \to \R{}$ is continuously differentiable and bounded below by $f_{\inf} \in \R{}$, and over which the constraint function $c : \R{n} \to \R{m}$ is continuously differentiable and bounded in $\ell_2$-norm in the sense that there exists $\kappa_c \in \R{}_{>0}$ such that $\|c(x)\|_2 \leq \kappa_c$ for all $x \in \Xcal$.  In addition, there exist constants $\kappa_{\nabla f} \in \R{}_{>0}$, $\kappa_J \in \R{}_{>0}$, $L_{\nabla f} \in \R{}_{>0}$, $\Lhat_J \in \R{}_{>0}$, and $\sigma_{\min} \in \R{}_{>0}$ such that one has that
  \begin{align*}
    \|\nabla f(x)\|_2 &\leq \kappa_{\nabla f} && \text{for all}\ \ x \in \Xcal, \\
    \|J(x)\|_2 &\leq \kappa_J && \text{for all}\ \ x \in \Xcal, \\
    \|\nabla f(x) - \nabla f(\xbar)\|_2 &\leq L_{\nabla f} \|x - \xbar\|_2 && \text{for all}\ \ (x,\xbar) \in \Xcal \times \Xcal, \\
    \|J(x) - J(\xbar)\|_2 &\leq \Lhat_J \|x - \xbar\|_2 && \text{for all}\ \ (x,\xbar) \in \Xcal \times \Xcal, \\
    \text{and}\ \ \sigma_1(J(x)) &\geq \sigma_{\min} && \text{for all}\ \ x \in \Xcal,
  \end{align*}
  where $\sigma_1(\cdot)$ yields the smallest singular value of its matrix argument.  Furthermore,
  \bequationNN
    \E[P_k g_k | \Fcal_k] = P_k \nabla f(x_k)\ \ \text{for all}\ \ k \in \N{},
  \eequationNN
  there exists a constant $M \in \R{}$ such that
  \bequationNN
    \E[\|P_k (g_k - \nabla f(x_k))\|_2^2 | \Fcal_k] \leq M^2\ \ \text{for all}\ \ k \in \N{},
  \eequationNN
  and there exist pairs of constants $(\rho_{\min},\rho_{\max}) \in (0,1] \times (0,1]$ and $(h_{\min},h_{\max}) \in \R{}_{>0} \times \R{}_{>0}$ such that $\rho_{\min} \leq \rho_k \leq \rho_{\max}$ and $h_{\min} \leq h_k \leq h_{\max}$ for all $k \in \N{}$.
\end{assumption}

It is well known (e.g., see \cite{BottCurtNoce2018}) that, under Assumption~\ref{ass.heavy-ball}, one has
\bequation\label{eq.f_Lipschitz}
  f(\xbar) - f(x) \leq \nabla f(x)^T(\xbar - x) + \thalf L_{\nabla f} \|\xbar - x\|_2^2\ \ \text{for all}\ \ (x,\xbar) \in \Xcal \times \Xcal,
\eequation
and, using norm inequalities, there exists $L_J \in \R{}_{>0}$ such that
\bequation\label{eq.J_Lipschitz}
  \|c(\xbar)\|_1 \leq \|c(x) + J(x)(\xbar - x)\|_1 + \thalf L_J \|\xbar - x\|_2^2\ \ \text{for all}\ \ (x,\xbar) \in \Xcal \times \Xcal.
\eequation

Next, let us show that an important Lipschitz continuity property holds.

\blemma\label{lem.projection_lipschitz}
  There exists $L_{P\nabla f} \in \R{}_{>0}$ such that for all $(x,\xbar) \in \Xcal \times \Xcal$ one has
  \bequationNN
    \|P(x) \nabla f(x) - P(\xbar) \nabla f(\xbar)\|_2 \leq L_{P\nabla f} \|x - \xbar\|_2.
  \eequationNN
\elemma
\proof{Proof.}
  Consider arbitrary $x \in \Xcal$.  Under Assumption~\ref{ass.heavy-ball}, the pseudoinverse of $J(x)$ is a right inverse and $J(x)^\dagger := J(x)^T (J(x) J(x)^T)^{-1}$, so $P(x) = I - J(x)^\dagger J(x)$.  Also, $P(x)$ is symmetric, and it is well known that $\|J(x)^\dagger\|_2 = \|J(x)^T (J(x) J(x)^T)^{-1}\|_2 \le \frac{1}{\sigma_{\min}}$ for any $x \in \mathcal{X}$.  Thus, for any $(x,\xbar) \in \Xcal \times \Xcal$, one has
  \begin{align*}
    \|P(x) - P(\xbar)\|_2
    &= \|P(x)(I-P(\xbar)) - (I - P(x))P(\xbar)\|_2 \\
    &= \|P(x)^T (I - P(\xbar))^T - (I-P(x)) P(\xbar)\|_2 \\
    &= \|P(x)^T J(\xbar)^T (J(\xbar)^\dagger)^T - J(x)^\dagger J(x) P(\xbar)\|_2 \\
    &= \|P(x)^T(J(\xbar) - J(x))^T (J(\xbar)^\dagger)^T - J(x)^\dagger (J(x) - J(\xbar)) P(\xbar) \|_2 \\
    &\leq \|P(x)\|_2 \|J(\xbar) - J(x)\|_2\|J(\xbar)^\dagger\|_2 + \|J(x)^\dagger\|_2 \|J(x) - J(\xbar)\|_2 \|P(\xbar)\|_2 \\
    &\leq (\|J(\xbar)^\dagger\|_2 + \|J(x)^\dagger\|_2) \|J(x) - J(\xbar)\|_2 \\
    &\leq \frac{2 \Lhat_J}{\sigma_{\min}} \|x - \xbar\|_2.
  \end{align*}
  Therefore, for any $(x,\xbar) \in \Xcal \times \Xcal$, one has
  \begin{align*}
    \|P(x)\nabla f(x) - P(\xbar)\nabla f (\xbar)\|_2
    =&\ \|P(x)\nabla f(x) - P(x)\nabla f(\xbar) + P(x)\nabla f(\xbar) - P(\xbar) \nabla f(\xbar)\|_2 \\
    \leq&\ \|P(x)\|_2 \|\nabla f(x) - \nabla f(\xbar)\|_2 + \|P(x) - P(\xbar)\|_2 \|\nabla f(\xbar)\|_2 \\
    \leq&\ \left (L_{\nabla f} +  \frac{2\Lhat_J \kappa_{\nabla f}}{\sigma_{\min}} \right) \|x - \xbar\|_2 =: L_{P\nabla f} \|x - \xbar\|_2,
    \end{align*}
  which completes the proof.
  \qed
\endproof

Next we state a couple of bounds pertaining to finite series.

\blemma\label{lem.betatoktimesk}
  Given $\beta \in (0,1)$ and $K \in \N{}$, one has that
  \begin{align*}
    \sum_{k=0}^{K-1} \beta^k = \frac{1-\beta^K}{1-\beta} \leq \frac{1}{1-\beta},\ \ \sum_{k=0}^{K-1} \beta^k k \leq \frac{\beta}{(1-\beta)^2},\ \ \text{and}\ \ \sum_{k=0}^{K-1} \beta^k k^2 \leq \frac{\beta(1+\beta)}{(1-\beta)^3}.
  \end{align*}
\elemma
\proof{Proof.}
  Each bound is straightforward to verify; e.g., for the second, see \cite[Lemma B.2]{defossezSimpleConvergenceProof2022}.
\endproof

We also state the following lemma, which follows easily from our prior observations under Assumption~\ref{ass.heavy-ball}.  Here and throughout the remainder of our analyses, we define $\phi : \R{n} \times \R{}_{>0} \to \R{}$ by $\phi(x,\tau) = \tau f(x) + \|c(x)\|_1$.

\blemma\label{lem.bounds}
  For all $k \in \N{}$, it follows that
  \begin{align} 
    \|v_k\|_2 &= \|\rho_k J_k^T (J_k J_k^T)^{-1} c_k\|_2 \leq \rho_k \kappa_c \sigma_{\min}^{-1}, \label{eq.v_bound} \\
    \E[\|u_k\|_2^2 | \Fcal_k] &= h_k^{-2} \E [\|P_k g_k\|_2^2 | \Fcal_k] \leq h_k^{-2} (\kappa_{\nabla f}^2 + M^2). \label{eq.u_bound}
  \end{align}
  In addition, for any $\tau \in \R{}_{>0}$, it follows for all $k \in \N{}$ that
  \bequation\label{eq.merit_reduction}
    \phi(x_k + \alpha d_k,\tau) - \phi(x_k,\tau) \leq \tau \alpha \nabla f(x_k)^T d_k + \|c_k + \alpha J_k d_k\|_1 - \|c_k\|_1 + \thalf \alpha^2 (\tau L_{\nabla f} + L_J) \|d_k\|_2^2.
  \eequation
\elemma
\proof{Proof.}
  Consider arbitrary $k \in \N{}$.  The bound \eqref{eq.v_bound} (respectively, \eqref{eq.u_bound}) follows from the definition of $v_k$ (respectively, $u_k$) in the algorithm and Assumption~\ref{ass.heavy-ball}.  In addition, under Assumption~\ref{ass.heavy-ball}, one has that
  \begin{align*}
    &\ \phi(x_k + \alpha d_k,\tau) - \phi(x_k,\tau) \\
    =&\ \tau (f(x_k + \alpha d_k) - f(x_k)) + \|c(x_k + \alpha d_k)\|_1 - \|c_k\|_1 \\ 
    \leq&\ \tau (\alpha \nabla f(x_k)^Td_k + \thalf \alpha^2 L_{\nabla f} \|d_k\|_2^2) + \|c_k + \alpha J_k d_k\|_1 + \thalf \alpha^2 L_J \|d_k\|_2^2 - \|c_k\|_1 \\
    =&\ \tau \alpha \nabla f(x_k)^T d_k + \|c_k + \alpha J_k d_k\|_1 - \|c_k\|_1 + \thalf \alpha^2 (\tau L_{\nabla f} + L_J) \|d_k\|_2^2, 
  \end{align*}
  which completes the proof.
  \qed
\endproof

For our next lemmas, first observe that for all $k \in \N{}$ one has that
\bequation\label{eq.r}
  r_k = u_k + \beta r_{k-1} = u_k + \beta (u_{k-1} + \beta r_{k-2}) = \sum_{i=0}^{k-1} \beta^{i} u_{k-i} = \sum_{i=1}^{k} \beta^{k-i} u_i.
\eequation

\blemma\label{lem.sqpmk.bound}
  For all $k \in \N{}$, it follows that
  \bequationNN
    \E[\|r_k\|_2^2] \leq \frac{\kappa_{\nabla f}^2 + M^2}{h_{\min}^2 (1-\beta)^2}.
  \eequationNN  
\elemma
\proof{Proof.}
  Consider arbitrary $k \in \N{}$.  By \eqref{eq.r}, Jensen's inequality, and Lemmas~\ref{lem.betatoktimesk} and \ref{lem.bounds} one finds that
  \begin{align*}
    \E[\|r_k\|_2^2] =
    \E \left[ \left\| \sum_{i=0}^{k-1} \beta^i u_{k-i} \right\|_2^2 \right]
    &=\ \E \left[ \( \sum_{i=0}^{k-1} \beta^i \)^2 \left\| \frac{\sum_{i=0}^{k-1} \beta^i u_{k-i}}{\sum_{i=0}^{k-1} \beta^i} \right\|_2^2 \right] \\
    &\leq\ \( \sum_{i=0}^{k-1} \beta^i \)^2 \frac{\sum_{i=0}^{k-1} \beta^i \E[\|u_{k-i}\|_2^2]}{\sum_{i=0}^{k-1} \beta^i} = \( \sum_{i=0}^{k-1} \beta^i \) \sum_{i=0}^{k-1} \beta^i \E[\|u_{k-i}\|_2^2] \\
    &\leq\ \frac{1}{1-\beta} \sum_{i=0}^{k-1} \beta^i h_k^{-2}(\kappa_{\nabla f}^2 + M^2) \leq \frac{\kappa_{\nabla f}^2 + M^2}{h_{\min}^2 (1-\beta)^2},
  \end{align*}
  as desired.
  \qed
\endproof

Next, we bound the expected value of the inner product between the true gradient of the objective function and the search direction in each iteration.

\blemma\label{lem.inner_product}
  For all $k \in \N{}$, it follows that
  \begin{align*}
    \E[\nabla f(x_k)^T d_k] \leq&\ - \sum_{i=0}^{k-1} \beta^i h_{\max}^{-1} \E[ \|P_{k-i} \nabla f(x_{k-i}) \|_2^2]  \\
    &\ + \frac{\rho_{\max} \kappa_{\nabla f} \E[\|c_k\|_2]}{\sigma_{\min}} +  \frac{\beta L_{P\nabla f} \alpha}{(1-\beta)^2} \( \frac{\rho_{\max}^2\kappa_c^2 (1-\beta)}{2\sigma_{\min}^2} + \frac{\kappa_{\nabla f}^2 + M^2}{ h_{\min}^2 (1-\beta)} \).
  \end{align*}
\elemma
\proof{Proof.}
  Consider arbitrary $k \in \N{}$.  By the definition of $d_k$, $P_k = P_k^T$, and \eqref{eq.r}, one finds that
  \begin{align}
    \nabla f(x_k)^T d_k
    =&\ \nabla f(x_k)^T (v_k + P_k r_k) \nonumber \\
    =&\ \nabla f(x_k)^T \(v_k +  P_k \sum_{i=0}^{k-1} \beta^i u_{k-i}\) \nonumber \\
    =&\ - \rho_k \nabla f(x_k)^T J_k^T (J_k J_k^T)^{-1} c_k + \sum_{i=0}^{k-1} \beta^i \nabla f(x_k)^T P_k u_{k-i} \nonumber \\
    =&\ - \rho_k \nabla f(x_k)^T J_k^T (J_k J_k^T)^{-1} c_k \nonumber \\
    &\ + \sum_{i=0}^{k-1} \beta^i \nabla f(x_{k-i})^T P_{k-i}^T u_{k-i} + \sum_{i=1}^{k-1} \beta^i (P_k \nabla f(x_k) - P_{k-i} \nabla f(x_{k-i}))^T u_{k-i}. \label{eq.descentdirection}
  \end{align}
  With respect to the first term in \eqref{eq.descentdirection}, it follows by Assumption~\ref{ass.heavy-ball} and Lemma~\ref{lem.bounds} (specifically, \eqref{eq.v_bound}) that
  \bequation\label{eq.firstterm}
    -\rho_k \nabla f(x_k)^T J_k^T (J_k J_k^T)^{-1} c_k \leq \frac{\rho_{\max} \kappa_{\nabla f} \|c_k\|_2}{\sigma_{\min}}.
  \eequation
  With respect to the second term in \eqref{eq.descentdirection}, it follows with Assumption~\ref{ass.heavy-ball} that for all $i \in \{0,\dots,k-1\}$ one has
  \begin{align}
    \E[\nabla f(x_{k-i})^T P_{k-i}^T u_{k-i} | \Fcal_{k-i}]
    &= \E[-h_{k-i}^{-1} \nabla f(x_{k-i})^T P_{k-i}^T P_{k-i} g_{k-i} | \Fcal_{k-i}] \nonumber \\
    &= -h_{k-i}^{-1} \nabla f(x_{k-i})^T P_{k-i}^T P_{k-i} \nabla f(x_{k-i}) \nonumber \\
    &= -h_{k-i}^{-1} \|P_{k-i} \nabla f(x_{k-i})\|_2^2 \leq -h_{\max}^{-1} \|P_{k-i} \nabla f(x_{k-i}) \|_2^2. \label{eq.secondterm}
  \end{align}
  With respect to the third term in \eqref{eq.descentdirection}, it follows by Lemma~\ref{lem.projection_lipschitz}, Jensen's inequality, Lemma~\ref{lem.bounds} (specifically, inequality~\eqref{eq.v_bound}), and $v_j$ and $P_jr_j$ being orthogonal for any $j \in \N{}$ that for all $i \in \{1,\dots,k-1\}$ one has
  \begin{align}
    \|P_k \nabla f(x_k) - P_{k-i} \nabla f(x_{k-i})\|_2^2
    &\leq L_{P\nabla f}^2 \|x_k - x_{k-i}\|_2^2 \nonumber \\
    &= L_{P\nabla f}^2 \left\| \alpha \sum_{j = k-i}^{k-1}(v_j + P_jr_j) \right\|_2^2 \nonumber \\ 
    &\leq L_{P\nabla f}^2 i \alpha^2 \(\sum_{j = k-i}^{k-1} \|v_j\|_2^2 + \sum_{j = k-i}^{k-1} \|P_jr_j\|_2^2 \) \nonumber \\
    &\leq L_{P\nabla f}^2 i \alpha^2 \(i \rho_{\max}^2 \kappa_c^2 \sigma_{\min}^{-2} + \sum_{j = k-i}^{k-1} \|P_jr_j\|_2^2 \). \label{eq.nablaf diff}
  \end{align}
  Consequently, along with the Cauchy-Schwarz inequality and since Young's inequality implies that $|ab| \leq \frac{\lambda}{2} |a|^2 + \frac{1}{2\lambda} |b|^2$ for any $(a,b) \in \R{} \times \R{}$ and $\lambda \in \R{}_{>0}$, one finds with $\lambda = \frac{1-\beta}{L_{P\nabla f} i \alpha}$ that
  \begin{align*}
    &\ \sum_{i=1}^{k-1} \beta^i (P_k \nabla f(x_k) - P_{k-i} \nabla f(x_{k-i}))^T u_{k-i} \\
    \leq&\ \sum_{i=1}^{k-1} \beta^i \|P_k \nabla f(x_k) - P_{k-i} \nabla f(x_{k-i})\|_2 \|u_{k-i}\|_2 \\
    \leq&\ \sum_{i=1}^{k-1} \beta^i \(\frac{1-\beta}{2 L_{P\nabla f} i \alpha} \|P_k \nabla f(x_k) - P_{k-i} \nabla f(x_{k-i})\|_2^2 + \frac{L_{P\nabla f} i \alpha}{2(1-\beta)} \|u_{k-i}\|_2^2 \) \\
    \leq&\ \sum_{i=1}^{k-1} \beta^i L_{P\nabla f} \alpha \( \frac{1-\beta}{2} \( i \rho_{\max}^2 \kappa_c^2 \sigma_{\min}^{-2} +  \sum_{j = k-i}^{k-1} \|P_jr_j\|_2^2\) + \frac{i}{2(1-\beta)} \|u_{k-i}\|_2^2 \).
  \end{align*}
  Taking total expectation and using Lemma~\ref{lem.betatoktimesk}, Lemma~\ref{lem.bounds} (specifically, \eqref{eq.u_bound}), and Lemma~\ref{lem.sqpmk.bound}, one finds that
  \begin{align}
    &\ \E\left[\sum_{i=1}^{k-1} \beta^i (P_k \nabla f(x_k) - P_{k-i} \nabla f(x_{k-i}))^T u_{k-i}\right] \nonumber \\
    \leq&\ \sum_{i=1}^{k-1} \beta^i L_{P\nabla f} \alpha \( \frac{1-\beta}{2} \(i \rho_{\max}^2 \kappa_c^2 \sigma_{\min}^{-2} + \sum_{j = k-i}^{k-1} \E\left[ \|P_jr_j\|_2^2 \right]\) + \frac{i}{2(1-\beta)} \E\left[\|u_{k-i}\|_2^2\right] \) \nonumber \\
    \leq&\ \sum_{i=1}^{k-1} \beta^i L_{P\nabla f} \alpha i \( \frac{1-\beta}{2} \( \rho_{\max}^2 \kappa_c^2 \sigma_{\min}^{-2} + \frac{ \kappa_{\nabla f}^2 + M^2}{h_{\min}^2(1-\beta)^2} \) + \frac{\kappa_{\nabla f}^2 + M^2}{2h_{\min}^2(1-\beta)} \) \nonumber \\
    =&\ L_{P\nabla f} \alpha \( \frac{ \rho_{\max}^2 \kappa_c^2 (1-\beta)}{2\sigma_{\min}^2} + \frac{ \kappa_{\nabla f}^2 + M^2}{ h_{\min}^2(1-\beta)}\) \sum_{i=1}^{k-1} \beta^i i \leq \frac{\beta L_{P\nabla f} \alpha}{(1-\beta)^2} \( \frac{\rho_{\max}^2 \kappa_c^2 (1-\beta)}{2\sigma_{\min}^2} + \frac{\kappa_{\nabla f}^2 + M^2}{ h_{\min}^2 (1-\beta)} \). \label{eq.thirdterm}
  \end{align}
  Taking expectation on both sides of \eqref{eq.descentdirection} and using \eqref{eq.firstterm}, \eqref{eq.secondterm}, and \eqref{eq.thirdterm}, the proof is complete.
\qed
\endproof

Now we bound the expected reduction in the merit function between two consecutive iterations.

\blemma\label{lem.merit_decrease}
  For any $\tau \in \R{}_{>0}$, it follows for all $k \in \N{}$ that
  \begin{align*}
    &\E[\phi(x_{k+1}) - \phi(x_k)] \\
    &\quad \leq - \alpha \( \tau \sum_{i=0}^{k-1} \beta^i h_{\max}^{-1} \E [ \|P_{k-i} \nabla f(x_{k-i}) \|_2^2 ] + \left(1 - \frac{\tau  \kappa_{\nabla f}\rho_{\max}}{\sigma_{\min}\rho_{\min}} \) \rho_{\min} \E[\|c_k\|_1] \) + \alpha^2 C, \\
    &\quad \text{where}\ \ C := \frac{\beta L_{P\nabla f} \tau}{1-\beta} \( \frac{\rho_{\max}^2 \kappa_c^2 }{2\sigma_{\min}^2} + \frac{\kappa_{\nabla f}^2 + M^2}{ h_{\min}^2 (1-\beta)^2}\) + \thalf (\tau L_{\nabla f} + L_J) \( \frac{\rho_{\max}^2 \kappa_c^2}{\sigma_{\min}^2} + \frac{\kappa_{\nabla f}^2 + M^2}{h_{\min}^2 (1-\beta)^2} \).
  \end{align*}
\elemma
\proof{Proof.}
  Consider arbitrary $\tau \in \R{}_{>0}$ and $k \in \N{}$.  By Lemma~\ref{lem.bounds}, $\alpha \in (0, 1]$, $\rho_k \in (0, 1]$, and
  \bequationNN
      J_k d_k = J_k v_k + J_k P_k r_k = J_k v_k = - \rho_k J_k J_k^T (J_k J_k^T)^{-1}c_k = -\rho_k c_k, 
  \eequationNN
  one has that
  \begin{align}
    \phi(x_{k+1}) - \phi(x_k)
      &\leq \tau \alpha \nabla f(x_k)^T d_k + \|c_k + \alpha J_k d_k\|_1 - \|c_k\|_1 + \thalf \alpha^2(\tau L_{\nabla f} + L_J) \|d_k\|_2^2 \nonumber \\
      &= \tau \alpha \nabla f(x_k)^T d_k - \alpha \rho_k \|c_k\|_1 + \thalf \alpha^2(\tau L_{\nabla f} + L_J) \|d_k\|_2^2. \label{eq.meritfuntiondiff1}
  \end{align}
  With respect to the last term, it follows from Lemma~\ref{lem.bounds}, Lemma~\ref{lem.sqpmk.bound}, and $v_k^TP_kr_k = 0$ that
  \begin{align*}
    \thalf \alpha^2(\tau L_{\nabla f} + L_J) \E[\|d_k\|_2^2]
    &=\ \thalf \alpha^2(\tau L_{\nabla f} + L_J) \E[\|v_k + P_kr_k\|_2^2] \\
    &\leq\ \thalf \alpha^2(\tau L_{\nabla f} + L_J) \( \frac{\rho_{\max}^2\kappa_c^2}{\sigma_{\min}^2} + \frac{ \kappa_{\nabla f}^2 + M^2}{h_{\min}^2 (1-\beta)^2} \).
  \end{align*}
  Now taking expectation on both sides of \eqref{eq.meritfuntiondiff1}, and using Lemma~\ref{lem.inner_product} and basic norm inequalities, one finds
  \begin{align*}
    &\ \E[\phi(x_{k+1}) - \phi(x_k)] \\
    \leq&\ - \tau \alpha \sum_{i=0}^{k-1} \beta^i h_{\max}^{-1} \E[\|P_{k-i} \nabla f(x_{k-i})\|_2^2] -\alpha \rho_{\min} \E[\|c_k\|_1] + \frac{\tau \alpha \rho_{\max} \kappa_{\nabla f} \E[\|c_k\|_1]}{\sigma_{\min}} \\
    &\ + \frac{\beta L_{P\nabla f} \alpha^2 \tau}{(1-\beta)^2} \( \frac{\rho_{\max}^2 \kappa_c^2 (1-\beta)}{2\sigma_{\min}^2} + \frac{\kappa_{\nabla f}^2 + M^2}{ h_{\min}^2 (1-\beta)}\) + \thalf \alpha^2(\tau L_{\nabla f} + L_J) \( \frac{\rho_{\max}^2 \kappa_c^2}{\sigma_{\min}^2} + \frac{\kappa_{\nabla f}^2 + M^2}{h_{\min}^2 (1-\beta)^2} \),
  \end{align*}
  which completes the proof.
  \qed
\endproof

\begin{theorem}\label{thm.SQP heavy-ball.convergence.rate}
  Suppose that Assumption~\ref{ass.heavy-ball} holds and define
  \begin{align} \label{eq.tau min}
    \tau := \frac{\sigma_{\min}\rho_{\min}}{\sigma_{\min}\rho_{\min} + \kappa_{\nabla f}\rho_{\max}} \implies 1 - \frac{\tau \kappa_{\nabla f}\rho_{\max}}{\sigma_{\min}\rho_{\min}} = \tau.
  \end{align}
  Then, for any $K \in \N{}$ and with $C \in \R{}_{>0}$ defined in Lemma~\ref{lem.merit_decrease}, it follows that
  \bequation\label{eq.convergence guarantee}
    \frac{1}{K} \sum_{k=1}^K \E[ h_{\max}^{-1} \|P_k \nabla f(x_k) \|_2^2 + \rho_{\min} \|c_k\|_1 ] \leq \frac{\phi(x_1) - \tau f_{\inf}}{\alpha \tau K} + \frac{\alpha C}{\tau}.
  \eequation
\end{theorem}
\proof{Proof.}
  Consider arbitrary $k \in \N{}$.  Given $\tau$ defined by \eqref{eq.tau min}, it follows from Lemma~\ref{lem.merit_decrease} that 
  \bequationNN
    \E[\phi(x_{k+1}) - \phi(x_k)] \leq - \alpha \tau \( \sum_{i=0}^{k-1} \beta^i h_{\max}^{-1} \E[ \|P_{k-i} \nabla f(x_{k-i}) \|_2^2 ] + \rho_{\min} \E[\|c_k\|_1] \) + \alpha^2 C.
  \eequationNN
  Summing this inequality over $k \in \{1, \dots, K\}$ and using $\phi(x) \geq \tau f_{\inf}$ for all $x \in \Xcal$ yields
  \bequation\label{eq.sumktoK}
    \frac{1}{K} \(\sum_{k=1}^K \sum_{i=0}^{k-1} \beta^i h_{\max}^{-1} \E[ \|P_{k-i} \nabla f(x_{k-i}) \|_2^2 ] + \sum_{k=1}^K \rho_{\min} \E[\|c_k\|_1] \) \leq \frac{\phi(x_1) - \tau f_{\inf}}{\alpha \tau K} + \frac{\alpha C}{\tau}. 
  \eequation
  With respect to the first term in the parentheses on the left-hand side, one finds that
  \begin{align*}
    \sum_{k=1}^K \sum_{i=0}^{k-1} \beta^i \E[ \|P_{k-i} \nabla f(x_{k-i}) \|_2^2 ]
    &= \sum_{k=1}^K \sum_{i=1}^k \beta^{k-i} \E[ \|P_i \nabla f(x_i) \|_2^2 ] \\
    &= \sum_{k=1}^K \E[ \|P_k \nabla f(x_k) \|_2^2 ]  \sum_{i=0}^{K-k} \beta^i \\
    &= \sum_{k=1}^K \E[ \|P_k \nabla f(x_k) \|_2^2 ]  \frac{1 - \beta^{K - k + 1}}{1-\beta} \\
    &\geq \sum_{k=1}^K \E[ \|P_k \nabla f(x_k)\|_2^2 ].
  \end{align*}
  This inequality, along with \eqref{eq.sumktoK}, completes the proof.
  \qed
\endproof

Theorem~\ref{thm.SQP heavy-ball.convergence.rate} shows that the average expected combination of the squared norm of the projected gradient and $\ell_1$-norm constraint violation over $K$ iterations decreases to $\alpha C/\tau$ at a rate of $\Ocal(1/K)$.  Indeed, the upper bound on the average expected combination can be made as small as desired: with $\{\rho_k\}$, $\{h_k\}$, and $\beta$ fixed, one can choose $\alpha$ small enough such that $\alpha C/\tau$ is as small as desired, then choose $K$ sufficiently large such that the first term on the right-hand side of \eqref{eq.convergence guarantee} is as small as desired.  In addition, $\rho_{\min}$ and $h_{\max}$ can be chosen to achieve any desired balance between $\E[\|P_k \nabla f(x_k)\|_2^2]$ and $\E[\|c_k\|_1]$ in the left-hand side of \eqref{eq.convergence guarantee}.

\subsection{Projected Stochastic Adam SQP}

Next, we consider an extension of the Adam method to the setting of equality-constrained optimization.  Specifically, we propose Algorithm~\ref{alg.padam_sqp}, where we note that for any equal-length real vectors $a$ and $b$ we use $a \circ b$ to denote their component-wise product, we use $e$ to denote a vector of all ones whose length is determined by the context in which it appears, and for any real vector~$v$ we use $\diag(v)$ to denote a diagonal matrix whose diagonal components are those of $v$ (in order).

We emphasize that, like the variant of Adam that is analyzed in \citep{defossezSimpleConvergenceProof2022}, Algorithm~\ref{alg.padam_sqp} involves a modified bias correction term in order to guarantee that a certain step size sequence (see the sequences $\{\eta_k\}$ and $\{\alpha_k\}$ below) is monotonically nondecreasing as $k \to \infty$.  As discussed in \citep{defossezSimpleConvergenceProof2022}, this variant for the unconstrained setting regularly yields comparable performance with the original Adam method.  We borrow this modified bias correction idea in order to model our analysis on that in \citep{defossezSimpleConvergenceProof2022}.

\begin{algorithm}[ht]
  \caption{Projected Stochastic Adam SQP}
  \label{alg.padam_sqp}
  \begin{algorithmic}[1]
    \Require $x_1 \in \R{n}$, $\{\rho_k\}$ and $\{h_k\}$ with $\rho_k \in (0,1]$ and $h_k \in \R{}_{>0}$ for all $k \in \N{}$, $\beta_1 \in [0,1)$, $\beta_2 \in (\beta_1, 1)$, $\alpha \in (0,1]$, and $\epsilon \in \R{}_{>0}$
    \State set $r_0 \gets 0 \in \R{n}$
    \State set $s_0 \gets 0 \in \R{n}$
    \For{all $k\in\N{}$}
      \State compute a stochastic gradient estimate $g_k \approx \nabla f(x_k)$
      \State compute $v_k \gets - \rho_k J_k^T \left(J_k J_k^T\right)^{-1} c_k$
      \State compute $u_k \gets - h_k^{-1}P_k g_k$, where $P_k:=I - J_k^T \left(J_k J_k^T\right)^{-1} J_k$
      \State set $r_k \gets \beta_1 r_{k-1} + u_k$ (first momentum)
      \State set $s_k \gets \beta_2 s_{k-1} + u_k \circ u_k$ (second momentum)
      \State set $\eta_k \gets \frac{(1-\beta_1)\sqrt{1-\beta_2^k}}{\sqrt{1-\beta_2}}$ (bias correction)
      \State set $d_k \gets v_k + \eta_k P_k \diag (s_k + \epsilon e)^{-1/2} r_k$
      \State Set $x_{k+1} \gets x_k + \alpha d_k$
    \EndFor
  \end{algorithmic}
\end{algorithm}

For our analysis of Algorithm~\ref{alg.padam_sqp}, we make the following assumption.

\begin{assumption}\label{ass.adam}
  With respect to Algorithm~\ref{alg.padam_sqp}, the conditions of Assumption~\ref{ass.heavy-ball} hold and, in addition,
  \bequationNN
    \|u_k\|_\infty \equiv \|h_k^{-1} P_k g_k\|_\infty \leq \sqrt{\frac{M^2 + \kappa_{\nabla f}^2}{h_k^2} - \epsilon}\ \ \text{for all}\ \ k \in \N{}. 
  \eequationNN
\end{assumption}

Before commencing our analysis, let us define a few quantities to simplify our expressions.  Let us also note that for any vector defined by the algorithm, we use a second subscript to denote a component index; e.g., for any $(k,i)$ we use $s_{k,i}$ to denote the $i$th component of the vector $s_k$.  Similarly, for the product $P_kg_k$, we use $[P_kg_k]_i$ to denote its $i$th component.  For the sake of simplicity, let us define the step size for $u_k$ as
\bequationNN
  \alpha_k := \alpha \eta_k = \frac{\alpha(1-\beta_1)\sqrt{1-\beta_2^k}}{\sqrt{1-\beta_2}}\ \ \text{for all}\ \ k \in \N{},
\eequationNN
and observe that it is nondecreasing as $k \to \infty$.  In addition, let us define for all $k \in \N{}$ the vectors $t_k$ and~$\ttilde_k$, where for all $i \in \{1,\dots,n\}$ the $i$th component of each vector is given by
\bequationNN
  t_{k,i} := \frac{r_{k,i}}{\sqrt{s_{k,i} + \epsilon}},\ \ \text{and}\ \ \ttilde_{k,i}:= \frac{u_{k,i}}{\sqrt{s_{k,i} + \epsilon}} =-\frac{[P_kg_k]_i}{h_k\sqrt{s_{k,i} + \epsilon}},
\eequationNN
respectively.  Let us also note that for all $(k,i)$ and $j \in \{1,\dots,k\}$ one has that
\bequationNN\label{eq.s_k_i}
    s_{k,i} = \sum_{l=1}^k \beta_2^{k-l} h_l^{-2}[P_lg_l]_i^2 = \beta_2^j s_{k-j,i} + \sum_{l=k-j+1}^k \beta_2^{k-l} h_l^{-2} [P_lg_l]_i^2,
\eequationNN
and for all $(k,i)$ and $j \in \{1,\dots,k\}$ let us define the related quantity
\bequation\label{eq.tilde_s_k_j_vector}
  \stilde_{k,j,i} = \beta_2^j s_{k-j,i} + \E \left[ \sum_{l=k-j+1}^k \beta_2^{k-l} h_l^{-2} [P_lg_l]_i^2 \Bigg| \Fcal_{k-{j+1}} \right].
\eequation
Also, similarly as for \eqref{eq.r} for the heavy-ball method, one finds here for all $k \in \N{}$ that
\bequation\label{eq.r_adam}
  r_k = u_k + \beta_1 r_{k-1} = u_k + \beta_1 (u_{k-1} + \beta_1 r_{k-2}) = \sum_{i=0}^{k-1} \beta_1^i u_{k-i} = \sum_{i=1}^k \beta_1^{k-i} u_i.
\eequation

We begin with two technical lemmas whose proofs can be found in \cite{defossezSimpleConvergenceProof2022}.

\blemma\label{lem.log}
  Let $(\beta_1,\beta_2,\epsilon)$ be given by Algorithm \ref{alg.padam_sqp} and let $\{a_k\}$ be any sequence of real numbers. For any $k \in \N{}$, with $b_k := \sum_{j=1}^k \beta_2^{k-j} a_j^2$ and $q_k := \sum_{j=1}^k \beta_1^{k-j} a_j$, one has that
  \begin{align*}
    \sum_{j=1}^k \frac{q_j^2}{b_j + \epsilon} &\leq \frac{1}{(1 - \beta_1) (1 - \frac{\beta_1}{\beta_2})} \( \log\(1 + \frac{b_k}{\epsilon}\) - k \log(\beta_2) \) \\
    \text{and}\ \ \sum_{j=1}^k \frac{a_j^2}{b_j + \epsilon} &\leq \log\(1 + \frac{b_k}{\epsilon}\) - k \log(\beta_2).
  \end{align*}
\elemma
\proof{Proof.}
  See \cite[Lemmas~5.2 and A.2]{defossezSimpleConvergenceProof2022}.
\endproof

\blemma\label{lem.series}
  For any $k \in \N{}$ and $\beta \in (0,1)$, it follows that
  \bequationNN
    \sum_{j=0}^{k-1} \beta^j \sqrt{j+1} \leq \frac{2}{(1 - \beta)^{3/2}}\ \ \text{and}\ \ \sum_{j=0}^{k-1} \beta^j \sqrt{j}(j+1) \leq \frac{4\beta}{(1 - \beta)^{5/2}}.
  \eequationNN
\elemma
\proof{Proof.}
  See \cite[Lemmas~A.3 and A.4]{defossezSimpleConvergenceProof2022}.
\endproof

Now, similarly as in Lemma~\ref{lem.inner_product} for the heavy-ball method in the previous subsection, we bound the expected inner product between the true gradient of the objective function and the search direction in each iteration.

\blemma \label{lem.inner_product_adam}
  For all $k \in \N{}$, it follows that
  \begin{align*}
    & \ \E[\nabla f(x_k)^T d_k] \\
    &\leq - \frac{\eta_k}{2h_{\max}} \sum_{j=0}^{k-1} \beta_1^j \E[ (P_{k-j}\nabla f(x_{k-j}))^T  \diag(\stilde_{k,j+1} + \epsilon e)^{-1/2}P_{k-j} \nabla f(x_{k-j}) ] \\
    & \  + \frac{\rho_{\max} \kappa_{\nabla f}\E[\|c_k\|_2]}{\sigma_{\min}} + \frac{3\eta_k \sqrt{M^2+\kappa_{\nabla f}^2} h_{\max}}{h_{\min}\sqrt{1-\beta_1}} \sum_{j=0}^{k-1} \left(\frac{\beta_1}{\beta_2} \right)^j\sqrt{j+1} \E[ \| \ttilde_{k-j}\|_2^2]\\
    & \  + \eta_k L_{P\nabla f}^2 \alpha^2 \rho_{\max}^2 \kappa_c^2 \sigma_{\min}^{-2} \frac{\sqrt{1-\beta_1}}{4\sqrt{M^2+\kappa_{\nabla f}^2}}\frac{\beta_1(1+\beta_1)}{(1-\beta_1)^3} \\
    &\ + \eta_k L_{P\nabla f}^2 \alpha_k^2 \frac{\sqrt{1-\beta_1}}{4\sqrt{M^2+\kappa_{\nabla f}^2}}\sum_{j=1}^{k-1} \E[\|t_{k-j}\|_2^2] \sum_{l=j}^{k-1} \beta_1^l \sqrt{l}.
  \end{align*}
\elemma
\proof{Proof.} Consider arbitrary $k \in \N{}$. By the definition of $d_k$, $P_k = P_k^T$, and \eqref{eq.r_adam}, one finds that
  \begin{align}
    \nabla f(x_k)^Td_k &= \nabla f(x_k)^T(v_k+\eta_k P_k  \diag(s_k + \epsilon e)^{-1/2}r_k) \nonumber \\
    &= \nabla f(x_k)^Tv_k - \eta_k \sum_{j=0}^{k-1} \beta_1^j h_{k-j}^{-1} \nabla f(x_k)^T P_k  \diag(s_k + \epsilon e)^{-1/2} P_{k-j} g_{k-j} \nonumber \\
    &= \underbrace{\nabla f(x_k)^Tv_k}_{A} - \underbrace{\eta_k \sum_{j=0}^{k-1} \beta_1^j h_{k-j}^{-1} (P_{k-j}\nabla f(x_{k-j}))^T  \diag(s_k + \epsilon e)^{-1/2} P_{k-j} g_{k-j}}_{B} \nonumber \\
    & \quad - \underbrace{\eta_k \sum_{j=0}^{k-1} \beta_1^j h_{k-j}^{-1} \left(P_k\nabla f(x_k) - P_{k-j}\nabla f(x_{k-j})\right)^T  \diag(s_k + \epsilon e)^{-1/2}P_{k-j} g_{k-j}}_{C}. \label{eq.lem91}
  \end{align}
  Term $A$ satisfies \eqref{eq.firstterm}. With respect to term $B$, one finds for any $j \in \{1,\dots,k-1\}$ that
  \begin{align*}
    (P_{k-j} \nabla f(x_{k-j}))^T \diag(s_k + \epsilon e)^{-1/2} P_{k-j} g_{k-j} =  \underbrace{(P_{k-j}\nabla f(x_{k-j}))^T \diag(\stilde_{k,j+1} + \epsilon e)^{-1/2}P_{k-j} g_{k-j}}_{B_1} \\
    + \underbrace{(P_{k-j}\nabla f(x_{k-j}))^T \left( \diag(s_k + \epsilon e)^{-1/2} - \diag(\stilde_{k,j+1} + \epsilon e)^{-1/2}\right)P_{k-j} g_{k-j}}_{B_2}.
  \end{align*}
  Thus, by the definition \eqref{eq.tilde_s_k_j_vector}, one finds for $B_1$ that
  \bequationNN
    \E[B_1] = \E[ \E[B_1 | \Fcal_{k-j}] ] = \E[ (P_{k-j}\nabla f(x_{k-j}))^T  \diag(\stilde_{k,j+1} + \epsilon e)^{-1/2}P_{k-j} \nabla f(x_{k-j}) ],
  \eequationNN
  and at the same time one finds for $B_2$ that, from \cite[Pages~19--20, Eq(A.27)]{defossezSimpleConvergenceProof2022} and Assumption~\ref{ass.adam}, one has
  \begin{align*}
    \E[|B_2|] \leq&\ \half \E[(P_{k-j} \nabla f(x_{k-j}))^T  \diag(\stilde_{k,j+1} + \epsilon e)^{-1/2} P_{k-j} \nabla f(x_{k-j}) ] \\
    &\ + \frac{2\sqrt{M^2+\kappa_{\nabla f}^2}}{h_{\min}\sqrt{1-\beta_1}\beta_2^j} \sqrt{j+1} \E[(P_{k-j}g_{k-j})^T  \diag(s_{k-j} + \epsilon e)^{-1}P_{k-j} g_{k-j}] \\
    =&\ \half \E[(P_{k-j}\nabla f(x_{k-j}))^T \diag(\stilde_{k,j+1} + \epsilon e)^{-1/2}P_{k-j} \nabla f(x_{k-j}) ] \\
    &\ + \frac{2\sqrt{M^2+\kappa_{\nabla f}^2}}{h_{\min}\sqrt{1-\beta_1}\beta_2^j} \sqrt{j+1} h_{k-j}^2\E[ \| \ttilde_{k-j}\|_2^2].
  \end{align*}
  From the definitions of $B_1$ and $B_2$, one finds that
  \begin{align*}
    B &= \eta_k \sum_{j=0}^{k-1} \beta_1^j h_{k-j}^{-1} (B_1 + B_2) \ge \eta_k \sum_{j=0}^{k-1} \beta_1^j h_{k-j}^{-1} (B_1 - |B_2|),
  \end{align*}
  so taking expectation and employing the above equation for $\E[B_1]$ and bound for $\E[|B_2|]$ yields
  \begin{align}
    \E[B]
    \geq&\ \eta_k \sum_{j=0}^{k-1} \beta_1^j h_{k-j}^{-1} (\E[B_1] - \E[|B_2|]) \nonumber \\ 
    \geq&\ \frac{\eta_k}{2h_{\max}} \sum_{j=0}^{k-1} \beta_1^j \E[ (P_{k-j}\nabla f(x_{k-j}))^T  \diag(\stilde_{k,j+1} + \epsilon e)^{-1/2}P_{k-j} \nabla f(x_{k-j}) ] \nonumber \\
    &\ -\frac{2\eta_k \sqrt{M^2+\kappa_{\nabla f}^2} h_{\max}}{h_{\min}\sqrt{1-\beta_1}} \sum_{j=0}^{k-1} \(\frac{\beta_1}{\beta_2} \)^j \sqrt{j+1} \E[\|\ttilde_{k-j}\|_2^2]. \label{eq.fec_B}
  \end{align}
  Now with respect to term $C$ in \eqref{eq.lem91}, with Lemma~\ref{lem.projection_lipschitz}, Jensen's inequality, Lemma~\ref{lem.bounds} (specifically, inequality~\eqref{eq.v_bound}), $v_{k-j}$ and $P_{k-j}t_{k-j}$ being orthogonal for any $k-l \in \N{}$, and the fact that $\{\alpha_k\}$ is nondecreasing, one finds (similarly as in \eqref{eq.nablaf diff}) that for any $j \in \{1,\dots,k-1\}$ one has
  \begin{align*}
    \|P_k\nabla f(x_k)-P_{k-j}\nabla f(x_{k-j})\|_2^2 &\leq L_{P\nabla f}^2\|x_k-x_{k-j}\|_2^2 \\
    &=  L_{P\nabla f}^2 \left\| \sum_{l=1}^j \left ( \alpha v_{k-l}+\alpha_{k-l} P_{k-l}t_{k-l} \right )\right\|_2^2 \\
    &\leq L_{P\nabla f}^2 j \left ( \sum_{l=1}^j \|\alpha v_{k-l} \|_2^2 + \sum_{l=1}^j \|\alpha_{k-l}P_{k-l}t_{k-l}\|_2^2 \right ) \\
    &\leq L_{P\nabla f}^2 j^2\alpha^2 \rho_{\max}^2 \kappa_c^2 \sigma_{\min}^{-2} + L_{P\nabla f}^2 j \alpha_{k-1}^2 \sum_{l=1}^j \|t_{k-l}\|_2^2 \\
    &\leq L_{P\nabla f}^2 j^2\alpha^2 \rho_{\max}^2 \kappa_c^2 \sigma_{\min}^{-2} + L_{P\nabla f}^2 j \alpha_k^2 \sum_{l=1}^j \|t_{k-l}\|_2^2.
  \end{align*}
  At the same time, for the vector $\diag(s_k + \epsilon e)^{-1/2}P_{k-j}g_{k-j}$, one finds for all $i \in \{1,\dots,n\}$ that
  \begin{align*} 
    \epsilon + s_{k,i} &= \epsilon + \beta_2^j s_{k-j,i}+\sum_{l=k-j+1}^k h_{l}^{-2} \beta_2^{k-l} (P_lg_l)_i^2 \geq \epsilon +\beta_2^j s_{k-j,i} \geq \beta_2^j(\epsilon+s_{k-j,i}) \\
    & \Longrightarrow \frac{(P_{k-j}g_{k-j})_i^2}{\epsilon + s_{k,i}} \leq \frac{1}{\beta_2^j}  \frac{(P_{k-j}g_{k-j})_i^2}{\epsilon + s_{k-j,i}}=\frac{h_{k-j}^2}{\beta_2^j}  \frac{(P_{k-j}g_{k-j})_i^2}{h_{k-j}^2(\epsilon + s_{k-j,i})} = \frac{h_{k-j}^2}{\beta_2^j} \ttilde^2_{k-j,i}.
  \end{align*}
  Thus, along with the Cauchy-Schwarz inequality and since Young's inequality implies that $|ab| \leq \frac{\lambda}{2} |a|^2 + \frac{1}{2\lambda} |b|^2$ for any $(a,b) \in \R{} \times \R{}$ and $\lambda \in \R{}_{>0}$, one finds with $\lambda = \frac{h_{\min}\sqrt{1-\beta_1}}{2\sqrt{M^2+\kappa_{\nabla f}^2}\sqrt{j+1}}$ that
  \begin{align*}
    &\ (P_k \nabla f(x_k) - P_{k-j} \nabla f(x_{k-j}))^T \diag(s_k + \epsilon e)^{-1/2} P_{k-j} g_{k-j} \\
    \leq &\ \|P_k \nabla f(x_k) - P_{k-j} \nabla f(x_{k-j})\|_2 \|\diag(s_k + \epsilon e)^{-1/2} P_{k-j} g_{k-j}\|_2 \\
    \leq &\ \frac{h_{\min}\sqrt{1-\beta_1}}{4\sqrt{M^2+\kappa_{\nabla f}^2}\sqrt{j+1}} \|P_k\nabla f(x_k) - P_{k-j} \nabla f(x_{k-j})\|_2^2 \\
    &\ + \frac{\sqrt{M^2+\kappa_{\nabla f}^2}\sqrt{j+1}}{h_{\min}\sqrt{1-\beta_1}} \|\diag(s_k + \epsilon e)^{-1/2} P_{k-j}g_{k-j}\|_2^2 \\
    \leq &\ \frac{h_{\min}\sqrt{1-\beta_1}}{4\sqrt{M^2+\kappa_{\nabla f}^2}\sqrt{j+1}} \left( L_{P\nabla f}^2 j^2\alpha^2 \rho_{\max}^2 \kappa_c^2 \sigma_{\min}^{-2} + L_{P\nabla f}^2 j \alpha_k^2 \sum_{l=1}^j \|t_{k-l}\|_2^2 \right) \\
    &\ + \frac{\sqrt{M^2+\kappa_{\nabla f}^2}\sqrt{j+1}}{h_{\min}\sqrt{1-\beta_1}} \frac{h_{k-j}^2}{\beta_2^j} \|\ttilde_{k-j}\|_2^2.
  \end{align*}
  Thus, using Lemma~\ref{lem.betatoktimesk}, one finds that
  \begin{align}
    &\ \E[-C] \nonumber \\
    \leq&\ \E[|C|] = \eta_k \sum_{j=0}^{k-1} \beta_1^j h_{k-j}^{-1} \E\left[ \left| (P_k \nabla f(x_k) - P_{k-j} \nabla f(x_{k-j}))^T \diag(s_k + \epsilon e)^{-1/2}P_{k-j} g_{k-j} \right| \right] \nonumber \\
    \leq&\ \eta_k L_{P\nabla f}^2 \alpha^2 \rho_{\max}^2 \kappa_c^2 \sigma_{\min}^{-2} \frac{\sqrt{1-\beta_1}}{4\sqrt{M^2+\kappa_{\nabla f}^2}}\sum_{j=0}^{k-1}  \frac{\beta_1^j j^2}{\sqrt{j+1}} \nonumber \\
    &\ + \eta_k L_{P\nabla f}^2 \alpha_k^2 \frac{\sqrt{1-\beta_1}}{4\sqrt{M^2+\kappa_{\nabla f}^2}}\sum_{j=0}^{k-1} \frac{\beta_1^j j}{\sqrt{j+1}}  \sum_{l=1}^j \E[\|t_{k-l}\|_2^2] \nonumber \\
    &\ + \frac{\eta_k \sqrt{M^2+\kappa_{\nabla f}^2}h_{\max}}{h_{\min}\sqrt{1-\beta_1}} \sum_{j=0}^{k-1} \left(\frac{\beta_1}{\beta_2}\right)^j \sqrt{j+1} \E[\|\ttilde_{k-j}\|_2^2] \nonumber \\
    \leq&\ \eta_k L_{P\nabla f}^2 \alpha^2 \rho_{\max}^2 \kappa_c^2 \sigma_{\min}^{-2} \frac{\sqrt{1-\beta_1}}{4\sqrt{M^2+\kappa_{\nabla f}^2}}\sum_{j=0}^{k-1} \beta_1^j j^2  + \eta_k L_{P\nabla f}^2 \alpha_k^2 \frac{\sqrt{1-\beta_1}}{4\sqrt{M^2+\kappa_{\nabla f}^2}}\sum_{j=0}^{k-1} \beta_1^j \sqrt{j} \sum_{l=1}^j \E[\|t_{k-l}\|_2^2] \nonumber \\
    &\ +  \frac{\eta_k\sqrt{M^2+\kappa_{\nabla f}^2}h_{\max}}{h_{\min}\sqrt{1-\beta_1}} \sum_{j=0}^{k-1} \left(\frac{\beta_1}{\beta_2}\right)^j \sqrt{j+1} \E[\|\ttilde_{k-j}\|_2^2] \nonumber \\
    \leq&\ \eta_k L_{P\nabla f}^2 \alpha^2 \rho_{\max}^2 \kappa_c^2 \sigma_{\min}^{-2} \frac{\sqrt{1-\beta_1}}{4\sqrt{M^2+\kappa_{\nabla f}^2}}\frac{\beta_1(1+\beta_1)}{(1-\beta_1)^3}  \nonumber \\
    &\ + \eta_k L_{P\nabla f}^2 \alpha_k^2 \frac{\sqrt{1-\beta_1}}{4\sqrt{M^2+\kappa_{\nabla f}^2}}\sum_{j=1}^{k-1} \E[\|t_{k-j}\|_2^2] \sum_{l=j}^{k-1} \beta_1^l \sqrt{l} \nonumber \\
    &\ + \frac{\eta_k \sqrt{M^2+\kappa_{\nabla f}^2}h_{\max}}{h_{\min}\sqrt{1-\beta_1}} \sum_{j=0}^{k-1} \left(\frac{\beta_1}{\beta_2}\right)^j \sqrt{j+1} \E[\|\ttilde_{k-j}\|_2^2]. \label{eq.fec_C}
  \end{align}
  Combining \eqref{eq.lem91} with the bounds \eqref{eq.firstterm}, \eqref{eq.fec_B}, and \eqref{eq.fec_C} completes the proof.
\qed
\endproof
\btheorem\label{thm.adam.convergence.rate}
  Suppose that Assumptions~\ref{ass.heavy-ball} and \ref{ass.adam} hold and define $\tau \in \R{}_{>0}$ by \eqref{eq.tau min} along with
  \begin{align*}
    G_1(\beta_1,\beta_2) &:= \frac{\tau \beta_1 L_{P\nabla f}^2}{(1-\beta_2)^{3/2}\sqrt{M^2+\kappa_{\nabla f}^2}},\ \ G_2(\beta_1,\beta_2) := \frac{(1-\beta_1) (\tau L_{\nabla f}+L_J)}{2(1-\beta_2)}, \\
    G_3(\beta_1,\beta_2) &:= \frac{\tau L_{P\nabla f}^2 \rho_{\max}^2 \kappa_c^2\sqrt{1-\beta_1}\beta_1(1+\beta_1)}{4\sqrt{M^2+\kappa_{\nabla f}^2}\sqrt{1-\beta_2}(1-\beta_1)^2\sigma_{\min}^{2}},\ \ \text{and}\ \ G_4(\beta_1,\beta_2) := \frac{(\tau L_{\nabla f}+L_J) \rho_{\max}^2 \kappa_c^2}{2\sigma_{\min}^2}.
  \end{align*}
  Then, for any $K \in \N{}$, it follows that
  \begin{align}
    &\ \frac{1}{K} \sum_{k=1}^K \( \frac{h_{\min}(1-\beta_1)}{2 h_{\max} \sqrt{M^2+\kappa_{\nabla f}^2}} \E[\|P_k\nabla f(x_k)\|_2^2] + \rho_{\min} \E[\|c_k\|_1] \) \nonumber \\
    \leq&\ \frac{\phi(x_1) - \tau f_{\inf}}{\alpha \tau K} + \frac{6h_{\max}\sqrt{M^2 + \kappa_{\nabla f}^2} \sqrt{1-\beta_1}} {h_{\min} \sqrt{1-\beta_2}} \frac{n}{(1-\frac{\beta_1}{\beta_2})^{3/2}} \( \frac{1}{K} \log \(1 + \frac{M^2 + \kappa_{\nabla f}^2}{h_{\min}^2 (1-\beta_2) \epsilon} \) -\log(\beta_2) \) \nonumber \\
    &\ + \(\frac{\alpha^2}{\tau} G_1(\beta_1,\beta_2) + \frac{\alpha}{\tau} G_2(\beta_1,\beta_2)\) \frac{n}{(1 - \frac{\beta_1}{\beta_2})} \( \frac{1}{K} \log \(1 + \frac{M^2 + \kappa_{\nabla f}^2}{h_{\min}^2 (1 - \beta_2) \epsilon} \) - \log(\beta_2) \) \nonumber \\
    &\ + \frac{\alpha^2}{\tau} G_3(\beta_1,\beta_2) + \frac{\alpha}{\tau} G_4(\beta_1,\beta_2). \label{eq.final_bound}
  \end{align}
\etheorem
\proof{Proof.}
Consider arbitrary $k \in \N{}$.  Similar to \eqref{eq.meritfuntiondiff1}, one has that 
  \begin{align}
    \phi(x_{k+1}) - \phi(x_k)
    &\leq \tau \alpha \nabla f(x_k)^Td_k - \alpha \rho_k \|c_k\|_1 + \thalf \alpha^2 (\tau L_{\nabla f} + L_J) \|d_k\|_2^2. \label{eq.adam_merit_decrease}
  \end{align}
  Let $\lambda_{\max}(\cdot)$ and $\lambda_{\min}(\cdot)$ denote the maximum and minimum eigenvalues, respectively, of a real symmetric matrix argument.  By Assumption~\ref{ass.adam} and Lemma~\ref{lem.betatoktimesk}, one finds for any $j \in \{1,\dots,k-1\}$ that
  \begin{align*}
    \lambda_{\max} \( \diag( \stilde_{k,j+1} + \epsilon e)^{1/2} \)
    &\leq \sqrt{\epsilon + \sum_{j=1}^k \beta_2^{k-j} \( \frac{M^2 + \kappa_{\nabla f}^2}{h_{\min}^2} - \epsilon \)} \leq \sqrt{\sum_{j=1}^k \beta_2^{k-j} \( \frac{M^2 + \kappa_{\nabla f}^2}{h_{\min}^2} \)} \\
    &= \frac{\sqrt{M^2 + \kappa_{\nabla f}^2}}{h_{\min}} \sqrt{\frac{1 - \beta_2^k}{1 - \beta_2}} = \frac{\alpha_k \sqrt{M^2 + \kappa_{\nabla f}^2}}{h_{\min} \alpha (1-\beta_1)}.
  \end{align*}
  Consequently, one finds that
  \begin{align}
    &\ \E[ (P_{k-j}\nabla f(x_{k-j}))^T \diag(\stilde_{k,j+1} + \epsilon e)^{-1/2} P_{k-j} \nabla f(x_{k-j}) ] \nonumber \\
    \geq&\ \E[ \lambda_{\min} ( \diag(\stilde_{k,j+1} + \epsilon e)^{-1/2} ) \|P_{k-j} \nabla f(x_{k-j})\|_2^2] \geq \frac{h_{\min} \alpha (1-\beta_1)}{\alpha_k \sqrt{M^2 + \kappa_{\nabla f}^2}} \E[ \|P_{k-j} \nabla f(x_{k-j})\|_2^2]. \label{eq.expectation_product_bound}
  \end{align}
  Taking total expectation~\eqref{eq.adam_merit_decrease}, using Lemmas~\ref{lem.bounds} and~\ref{lem.inner_product_adam} and \eqref{eq.expectation_product_bound}, one finds that
  \begin{align}
    &\ \E[\phi(x_{k+1})-\phi(x_k)] \nonumber \\
    \leq&\ \tau \alpha \E[\nabla f(x_k)^Td_k] - \alpha \rho_k \E[\|c_k\|_1] + \thalf \alpha^2 (\tau L_{\nabla f}+L_J) \E[\|d_k\|_2^2] \nonumber \\
    \leq&\ \tau \alpha \E[\nabla f(x_k)^Td_k] - \alpha \rho_k \E[\|c_k\|_1] + \thalf \alpha^2 (\tau L_{\nabla f}+L_J) ( \E[\|v_k\|_2^2] + \eta_k^2\E[\|t_k\|_2^2]) \nonumber \\
    \leq&\ - \frac{\tau h_{\min} \alpha (1-\beta_1)}{2h_{\max} \sqrt{M^2+\kappa_{\nabla f}^2}} \sum_{j=0}^{k-1} \beta_1^j \E[\| P_{k-j}\nabla f(x_{k-j}) \|_2^2] + \tau \alpha \frac{\rho_{\max} \kappa_{\nabla f} \E[\|c_k\|_1]}{\sigma_{\min}} \nonumber \\
    &\ + \tau \alpha_k \frac{3h_{\max} \sqrt{M^2 + \kappa_{\nabla f}^2}}{h_{\min} \sqrt{1 - \beta_1}}  \sum_{j=0}^{k-1} \(\frac{\beta_1}{\beta_2} \)^j \sqrt{j+1} \E[ \|\ttilde_{k-j}\|_2^2] \nonumber \\
    &\ + \tau \alpha_k L_{P\nabla f}^2 \alpha^2 \rho_{\max}^2 \kappa_c^2 \sigma_{\min}^{-2} \frac{\sqrt{1-\beta_1}}{4\sqrt{M^2 + \kappa_{\nabla f}^2}} \frac{\beta_1(1+\beta_1)}{(1-\beta_1)^3} \nonumber \\
    &\ + \tau \alpha_k^3 L_{P\nabla f}^2 \frac{\sqrt{1-\beta_1}}{4\sqrt{M^2 + \kappa_{\nabla f}^2}} \sum_{j=1}^{k-1} \E[\|t_{k-j}\|_2^2] \sum_{l=j}^{k-1} \beta_1^l \sqrt{l} \nonumber \\
    &\ - \alpha \rho_{\min} \E[\|c_k\|_1] + \thalf \alpha^2 (\tau L_{\nabla f}+L_J) \(\rho_{\max}^2 \kappa_c^2 \sigma_{\min}^{-2} + \eta_k^2 \E[\|t_k\|_2^2] \). \nonumber 
  \end{align}
  Summing over $k \in \{1,\dots,K\}$ and using the fact that $\{\alpha_k\}$ is nondecreasing, one has
  \begin{align}
    &\ \underbrace{\frac{\tau h_{\min} \alpha (1-\beta_1)}{2 h_{\max} \sqrt{M^2 + \kappa_{\nabla f}^2}} \sum_{k=1}^K \sum_{j=0}^{k-1} \beta_1^j \E[ \| P_{k-j}\nabla f(x_{k-j}) \|_2^2]}_{A} + \alpha \rho_{\min} \(1 - \frac{\tau \kappa_{\nabla f}\rho_{\max}}{\sigma_{\min}\rho_{\min}}\) \sum_{k=1}^K \E[\|c_k\|_1] \nonumber \\
    \leq&\ \phi(x_1) - \tau f_{\inf} + \underbrace{\frac{3h_{\max}\sqrt{M^2 + \kappa_{\nabla f}^2} \tau \alpha_K} {h_{\min} \sqrt{1-\beta_1}} \E \left[ \sum_{k=1}^K \( \sum_{j=0}^{k-1} \( \frac{\beta_1}{\beta_2} \)^j \sqrt{j+1} \|\ttilde_{k-j}\|_2^2 \) \right]}_{B} \nonumber \\
    &\ + \underbrace{\frac{\tau \alpha_K^3 L_{P\nabla f}^2 \sqrt{1-\beta_1}}{4 \sqrt{M^2 + \kappa_{\nabla f}^2}} \E \left[ \sum_{k=1}^K \( \sum_{j=1}^{k-1} \|t_{k-j}\|_2^2 \sum_{l=j}^{k-1} \beta_1^l \sqrt{l} \) \right]}_{C} + \frac{\alpha_K^2 (\tau L_{\nabla f} + L_J)}{2} \E \left[ \sum_{k=1}^K \|t_k\|_2^2 \right] \nonumber \\
    &\ + K \( \tau\alpha_K L_{P\nabla f}^2 \alpha^2 \rho_{\max}^2 \kappa_c^2 \sigma_{\min}^{-2} \frac{\sqrt{1-\beta_1}}{4\sqrt{M^2 + \kappa_{\nabla f}^2}} \frac{\beta_1(1+\beta_1)}{(1-\beta_1)^3} + \frac{\alpha^2 (\tau L_{\nabla f} + L_J) \rho_{\max}^2 \kappa_c^2}{2\sigma_{\min}^2} \). \label{eq.ABCDEG}
  \end{align}
  With respect to the term $A$, one finds by Lemma~\ref{lem.betatoktimesk} that
  \begin{align*}
    A
    &= \frac{\tau h_{\min} \alpha (1-\beta_1)}{2h_{\max} \sqrt{M^2 + \kappa_{\nabla f}^2}} \sum_{k=1}^K \sum_{j=0}^{k-1} \beta_1^j \E \left[ \| P_{k-j} \nabla f(x_{k-j}) \|_2^2 \right] \\
    &= \frac{\alpha h_{\min} \tau}{2 h_{\max} \sqrt{M^2 + \kappa_{\nabla f}^2}} \sum_{k=1}^K (1-\beta_1^{K-k+1}) \E[\|P_k\nabla f(x_k)\|_2^2].
  \end{align*}
  With respect to the term $B$, one finds with Lemma \ref{lem.series} and $\beta_1 < \beta_2$ that 
  \begin{align*}
    B
    &= \frac{3h_{\max} \sqrt{M^2 + \kappa_{\nabla f}^2} \tau \alpha_K} {h_{\min} \sqrt{1-\beta_1}} \E \left[ \sum_{k=1}^K \|\ttilde_k\|_2^2 \sum_{j=k}^K \( \frac{\beta_1}{\beta_2} \)^{j-k} \sqrt{1+j-k} \right] \\
    &\leq \frac{6 h_{\max} \sqrt{M^2 + \kappa_{\nabla f}^2} \tau \alpha_K} {h_{\min} \sqrt{1-\beta_1}} \frac{1}{(1-\frac{\beta_1}{\beta_2})^{3/2}} \E \left[ \sum_{k=1}^K \|\ttilde_k\|_2^2 \right].
  \end{align*}
  With respect to the term $C$, one finds with Lemma \ref{lem.series} that
  \bequationNN
    C \leq \frac{\tau \alpha_K^3 L_{P\nabla f}^2 \sqrt{1-\beta_1}}{4 \sqrt{M^2 + \kappa_{\nabla f}^2}} \E \left[ \sum_{k=1}^K \|t_k\|_2^2 \sum_{l=0}^{K-1} \beta_1^l \sqrt{l}(l+1) \right] \leq \frac{\tau \alpha_K^3 L_{P\nabla f}^2}{\sqrt{M^2 + \kappa_{\nabla f}^2}} \frac{\beta_1}{(1-\beta_1)^2} \E \left[ \sum_{k=1}^K \|t_k\|_2^2 \right].
  \eequationNN
  Now it follows from Assumption~\ref{ass.adam} and Lemma~\ref{lem.betatoktimesk} that $s_{K,i} \leq \frac{M^2+\kappa_{\nabla f}^2}{h_{\min}^{2}(1-\beta_2)}$. Thus, with Lemma~\ref{lem.log},
  \begin{align*}
    \sum_{k=1}^K \|t_k\|_2^2
    &= \sum_{i=1}^n \sum_{k=1}^K t_{k,i}^2 = \sum_{i=1}^n \sum_{k=1}^K\frac{r_{k,i}^2}{s_{k,i} + \epsilon} \\
    &\leq \sum_{i=1}^n \frac{1}{(1-\beta_1)(1-\frac{\beta_1}{\beta_2})} \( \log \( 1 + \frac{s_{K,i}}{\epsilon} \) - K \log(\beta_2) \) \\
    &\leq \frac{n}{(1-\beta_1)(1-\frac{\beta_1}{\beta_2})} \( \log \( 1 + \frac{M^2 + \kappa_{\nabla f}^2}{h_{\min}^{2} (1-\beta_2) \epsilon} \) - K \log(\beta_2) \)
  \end{align*}
   and
  \begin{align*}
    \sum_{k=1}^K \|\ttilde_k\|_2^2 = \sum_{i=1}^n \sum_{k=1}^K \ttilde_{k,i}^2 = \sum_{i=1}^n \sum_{k=1}^K \frac{u_{k,i}^2}{ s_{k,i} + \epsilon }
    &\leq \sum_{i=1}^n \( \log \( 1 + \frac{s_{K,i}}{\epsilon} \) - K \log(\beta_2) \) \\
    &\leq n \( \log \(1 + \frac{M^2 + \kappa_{\nabla f}^2}{h_{\min}^{2}(1-\beta_2)\epsilon} \) - K \log(\beta_2) \).
  \end{align*}
  Hence, with $\alpha_K \leq \alpha \frac{1-\beta_1}{\sqrt{1-\beta_2}}$ and the above bounds for $A$, $B$, and $C$, one has from \eqref{eq.ABCDEG} that
  \begin{align*}
    &\ \frac{\alpha \tau h_{\min}}{2 h_{\max} \sqrt{M^2 + \kappa_{\nabla f}^2}} \sum_{k=1}^K (1-\beta_1^{K-k+1}) \E[\|P_k\nabla f(x_k)\|_2^2] + \alpha \rho_{\min} \(1 - \frac{\tau \kappa_{\nabla f}\rho_{\max}}{\sigma_{\min}\rho_{\min}}\) \sum_{k=1}^K \E[\|c_k\|_1] \\
    \leq&\ \phi(x_1) - \tau f_{\inf} + \frac{6 h_{\max} \sqrt{M^2 + \kappa_{\nabla f}^2} \tau \alpha \sqrt{1-\beta_1}} {h_{\min} \sqrt{1-\beta_2}} \frac{n}{(1 - \frac{\beta_1}{\beta_2})^{3/2}} \( \log \(1 + \frac{M^2 + \kappa_{\nabla f}^2}{h_{\min}^{2} (1-\beta_2) \epsilon} \) -K \log(\beta_2) \) \\
    &\ + \( \alpha^3 G_1(\beta_1,\beta_2) + \alpha^2 G_2(\beta_1,\beta_2) \) \frac{n}{(1-\frac{\beta_1}{\beta_2})} \( \log \(1 + \frac{M^2 + \kappa_{\nabla f}^2}{h_{\min}^2 (1-\beta_2) \epsilon} \) - K \log(\beta_2) \) \\
    &\ + K \( \alpha^3 G_3(\beta_1,\beta_2) + \alpha^2 G_4(\beta_1,\beta_2) \).
  \end{align*}
  Diving both sides by $\alpha \tau K$, and using \eqref{eq.tau min} along with $(1-\beta_1^{K-k+1}) \geq (1-\beta_1)$, yields the result.
  \qed
\endproof

This theorem shows that the long-run average of a positive combination of $\E[\|P_k\nabla f(x_k)\|_2^2]$ and $\E[\|c_k\|_1]$ can be made as small as desired. Consider the right-hand side of \eqref{eq.final_bound}, which can be written as the sum of
\begin{align*}
  A &:= \frac{\phi(x_1) - \tau f_{\inf}}{\alpha \tau K}, \\
  B &:= \frac{6h_{\max}\sqrt{M^2 + \kappa_{\nabla f}^2} \sqrt{1-\beta_1}n\log \(1 + \frac{M^2 + \kappa_{\nabla f}^2}{h_{\min}^2 (1-\beta_2) \epsilon} \)} {h_{\min} \sqrt{1-\beta_2}(1-\frac{\beta_1}{\beta_2})^{3/2}K}, \\
  C &:= \frac{6h_{\max}\sqrt{M^2 + \kappa_{\nabla f}^2} \sqrt{1-\beta_1}n\( -\log(\beta_2) \) } {h_{\min} \sqrt{1-\beta_2}(1-\frac{\beta_1}{\beta_2})^{3/2}}, \\
  D &:= \(\frac{\alpha^2}{\tau} G_1(\beta_1,\beta_2) + \frac{\alpha}{\tau} G_2(\beta_1,\beta_2)\) \frac{n}{(1 - \frac{\beta_1}{\beta_2})} \( \frac{1}{K} \log \(1 + \frac{M^2 + \kappa_{\nabla f}^2}{h_{\min}^2 (1 - \beta_2) \epsilon} \) - \log(\beta_2) \), \\
  \text{and}\ \ E &:= \frac{\alpha^2}{\tau} G_3(\beta_1,\beta_2) + \frac{\alpha}{\tau} G_4(\beta_1,\beta_2).
\end{align*}
Supposing that $\{\rho_k\}$, $\{h_k\}$, $\beta_1$, and $\epsilon$ are set and fixed, one can choose the remaining inputs $\beta_2$, $\alpha$, and $K$ to make each of the above terms as small as desired. Specifically, consider first the term $C$. Observe that $\frac{-\log(\beta_2)}{\sqrt{1-\beta_2}}$ and $(1-\frac{\beta_1}{\beta_2})^{-3/2}$ are both monotonically decreasing in $\beta_2 \in (\beta_1, 1)$ as $\beta_2 \to 1$ and the former decreases to $0$ and the latter to $(1-\beta_1)^{-3/2}$, and hence $C$ is monotonically decreasing in $\beta_2$ over this range and decreases to $0$ as $\beta_2 \to 1$. Thus, one can first choose $\beta_2$ close enough to $1$ such that $C$ is as small as desired, and then consider the value of $\beta_2$ as fixed. Next, one can choose $\alpha \in (0, 1]$ small enough such that $D$ and $E$ are as small as desired, and then consider the value of $\alpha$ as fixed. Note that $D$ could be further reduced by increasing $K$. The last parameter to choose is $K$, where one can choose $K$ large enough such that $A$ and $B$ are as small as desired.  In summary, by first choosing $\beta_2$ close enough to $1$, then choosing $\alpha$ small enough, and finally choosing $K$ large enough, one can make the right-hand side of \eqref{eq.final_bound} as small as desired. Consequently, the long-run average of the linear combination of $\E[\|P_k\nabla f(x_k)\|_2^2]$ and $\E[\|c_k\|_1]$ can be made as small as desired.

\section{Informed Supervised Machine Learning and Implementation Details}\label{sec.implementation}

Let us now discuss a methodology for which the algorithms proposed and analyzed in the previous section are particularly well suited. The methodology involves incorporating prior knowledge into a supervised learning process through \emph{hard constraints} that are imposed \emph{during training only}.  Both of these highlighted aspects are critical for its effectiveness.  The methodology's use of \emph{hard constraints} is in contrast to previously proposed methodologies that incorporate prior knowledge through either (a) soft constraints \citep{Zhu2019} (i.e., through regularization/penalty terms in the objective function) or (b) designing the prediction function to incorporate knowledge directly \citep{Chalapathi2024,Negiar2023} (e.g., through neural network layers for which a forward pass requires solving a set of equations or even an optimization problem).  By imposing such constraints \emph{during training only}, one can avoid having the trained network require expensive operations for each forward pass.  Another key feature of this methodology is that one does not solve the hard-constrained training problem with a penalty-based (e.g., augmented Lagrangian) method.  This feature is also critical for the effectiveness of the methodology.

The supervised training of a machine learning model involves solving an optimization problem over a set of parameters of a prediction function, call it $p : \R{n_f} \times \R{n} \to \R{n_o}$, where $n_f$ is the number of features in an input,~$n$ is the dimension of the training/optimization problem, and $n_o$ is the dimension of the output.  Denoting known input-output pairs in the form $(a,b) \in \R{n_f} \times \R{n_o}$ and given a loss function $\ell : \R{n_o} \times \R{n_o} \to \R{}$, the training/optimization problem can be viewed in expected-loss or empirical-loss minimization form, i.e.,
\bequationNN
  \min_{x \in \R{n}}\ \textstyle \int_{\Acal \times \Bcal} \ell(p(a,x),b) \text{d}\P_{A,B}(a,b)\ \ \approx\ \ \min_{x \in \R{n}}\ \frac{1}{N} \sum_{i=1}^N \ell(p(a_i,x),b_i),
\eequationNN
where $\Acal$ is the input domain, $\Bcal$ is the output domain, $\P_{A,B}$ is the input-output probability function, and $\{(a_i,b_i)\}_{i=1}^N \subset \R{n_f} \times \R{n_o}$.  In our setting, the problem has (hard) constraints on $x$ as well.  Generally, these can be formulated in various ways; e.g., expectation, probabilistic, or almost-sure constraints.  We contend that for many informed-learning problems---such as for many physics-informed learning problems, as we discuss below---a fixed, small number of constraints suffices to improve training.  Given a (small) number $m$ of input-output pairs $\{(a_i^c,b_i^c)\}_{i=1}^{m}$, the constraints may take the form
\bequationNN
  \phi_i(p(a_i^c,x),b_i^c,\dots) = 0\ \ \text{for all}\ \ i \in \{1,\dots,m\},
\eequationNN
where the arguments to the constraint functions $\{\phi_i\}$ may include additional terms, such as derivatives of the prediction function with respect to inputs and/or model weights; see~\S\ref{sec.experiments} for specific examples.

Our algorithms from the previous section can be employed in numerous informed-learning contexts (e.g., fair learning \cite{CurtLiuRobi23,MDonini_etal_2018,JKomiyama_etal_2018,MBZafar_etal_2017,MBZafar_etal_2017b,ZafaValeGomeGumm19}).  For this work, we tested our approach on a few physics-informed learning problems.  We emphasize that our goal here is not to test huge-scale, state-of-the-art techniques for physics-informed learning.  Rather, we take a few physics-informed learning test problems and train relatively straightforward neural networks in order to demonstrate the relative performance of our proposed algorithms with a soft-constrained approach with Adam scaling~\citep{KingBa14} and the hard-constrained approach with Adam scaling from~\citep{Marquez-Neila}.  The relative performance of the algorithms would be similar if we were to train much more sophisticated and large-scale neural networks that are being developed in state-of-the-art physics-informed learning.  For more on physics-informed learning we direct the reader to, e.g., \citep{Cuomo2022,Karniadakis2021,lagaris1998artificial,Raissi2019,Takamoto2023,Wang2021,Wang2023}. The work \citep{chen2024physics} lies in the physics-informed learning with hard constraints, but is restricted to the hard constraints that the PDE inputs and solutions are linearly related, whereas our method handles general nonlinear constraints. They enforce feasibility via projection, while we allow infeasible iterates, using projection only for momentum. Thus, we do not compare our method with theirs.

Let us now provide an overview of the setting of physics-informed learning that we consider in our experiments.  A parametric partial differential equation (PDE) can be written generically as $\Fcal(\phi, u) = 0$, where $(\Phi, \Ucal, \Vcal)$ is a triplet of Banach spaces, $\Fcal: \Phi \times \Ucal \to \Vcal$ is a differential operator, $\phi \in \Phi$ represents PDE parameters, and $u \in \Ucal$ denotes a solution of the PDE corresponding to $\phi$.  The aim is to train a model to learn a mapping from the PDE parameters to a corresponding solution.  Let such a mapping be denoted as $\Gcal : \Phi \times \R{d} \times \R{n} \to \Ucal$, the inputs to which are PDE parameters, a vector encoding information about the domain of the PDE solution about which one aims to make a prediction (e.g., temporal and/or spatial coordinates), and, say, neural-network model parameters, and the output is a predicted solution value.

For training a model to solve the PDE with potentially no known solution values (see \citep{Karniadakis2021}), one can consider a set of training inputs $\{(\phi_i, y_i)\}_{i \in S_1}$ and minimize the average PDE residual over the training inputs.  Assuming that, in addition, one has access to observed and/or computed solution data in the form of tuples $\{(\phi_i, y_i, u_i)\}_{i \in S_2}$, one can also aim to minimize the differences between known and predicted solution values.  Mathematically, these aims can be expressed as finding $x$ to minimize
\bequation\label{eq.const}
  \frac{1}{|S_1|} \sum_{i \in S_1} \|\Fcal(\phi_i, \Gcal(\phi_i,y_i,x))\|_2^2\ \ \text{and/or}\ \ \frac{1}{|S_2|} \sum_{i \in S_2} \|u_i - \Gcal(\phi_i,y_i,x)\|_2^2.
\eequation
Note that the $\phi_i$ and/or $y_i$ elements in $\{(\phi_i, y_i)\}_{i \in S_1}$ may be the same or different from those in $\{(\phi_i, y_i, u_i)\}_{i \in S_2}$.  Additional terms may also be used for training, e.g., pertaining to initial and/or boundary conditions, or pertaining to \emph{partial} physics information.  For example, in \S\ref{sec.chemistry}, we train a model for which it is known that a mass-balance equation should hold, so our training problem involves residuals for the known mass-balance equation, even though this only defines the physics partially.  Overall, if one combines all learning aims into a single objective function---say, with a linear combination involving weights for the different objective terms---then one is employing a \emph{soft-constrained} approach.  We contend that a more effective approach can be to take at least a subset of terms and impose them as \emph{hard constraints} during training.  For example, with respect to the aims in (\ref{eq.const}), one might impose constraints such as $\Fcal(\phi_i, \Gcal(\phi_i,y_i,x)) = 0$ for some $i \in S_1$ and/or $u_i = \Gcal(\phi_i,y_i,x)$ for some $i \in S_2$.  Our experiments show the benefits of this idea.

\section{Numerical Experiments}\label{sec.experiments}

In this section, we present the results of numerical experiments that compare the performance of our proposed algorithms (\texttt{SQP-Heavyball(con)} and \texttt{SQP-Adam(con)}, respectively, where (``con'' stands for ``constrained'') versus a soft-constrained approach with Adam scaling (\texttt{Adam(unc)} for ``unconstrained'') \citep{Karniadakis2021} and a hard-constrained approach with \emph{projection-less} Adam scaling (\texttt{Adam(con)}) \citep{Marquez-Neila}.  We consider four test problems.  A few of them---namely, our 1D spring, 1D Burgers' equation, and 2D Darcy flow problems---have been seen in the literature; see \citep{li2021fourier, Negiar2023}.  We also consider a problem from chemical engineering, a modified version of a reaction network proposed in \citep{Gupta2016}.  To ensure a fair comparison, for each test problem, \texttt{SQP-Heavyball(con)}, \texttt{SQP-Adam(con)}, \texttt{Adam(unc)}, and \texttt{Adam(con)} use the same objective function, which is usually the data-fitting loss plus PDE residual loss, while the hard-constrained methods \texttt{SQP-Heavyball(con)}, \texttt{SQP-Adam(con)}, and \texttt{Adam(con)} impose additional constraints: the PDE residuals are zero at some input data points.  Further details are provided in each problem's subsection.  
The software uses PyTorch (BSD-3 license).  For all experiments, the parameters $\beta = 0.9$, $\beta_1=0.9$, $\beta_2 = 0.999$, and $\epsilon = 10^{-7}$ were used; see Algorithms~\ref{alg.heavy-ball} and~\ref{alg.padam_sqp}.  Our numerical experiments were performed using Google CoLaboratory\textsuperscript{TM} L4 GPU platforms.

\subsection{1D Spring} \label{subsec: 1dspring}

Our first test problem aims to predict the movement of a damped harmonic (mass-spring) oscillator \citep{harmonic.oscillator} under the influence of a restoring force and friction.  For simplicity, our aim was to train a model to predict the movement for known parameters and a single initial condition.  (Our later test problems involved more complicated situations; this simple problem and the case of only a single initial condition serves as a good starting point for comparison.)  The spring can be described by a linear, homogeneous, second-order ordinary differential equation with constant coefficients, namely, $m \frac{d^{2}u(t)}{dt^{2}} + \mu \frac{du(t)}{dt} + ku(t) = 0$ over $t \in [0,1]$, where we fixed the mass $m=1$, friction coefficient $\mu=4$, and spring constant $k=400$.  This corresponds to an under-damped state for which the exact solution with amplitude $A$ and phase $\phi$ is well known to be $u(t)= e^{-\delta t}(2A\cos(\phi + t\sqrt{w_{0}^{2} - \delta^{2}}))$, where $\delta = \mu/(2m)$ and $w_0 = \sqrt{k/m}$.

Our aim was to train a neural network with the known ODE and a few observed solution values to be able to predict the height of the spring at any time $t \in [0,1]$.  We used a fully connected neural network with~1 input neuron (corresponding to $t$), 3 hidden layers with 32 neurons each, and~1 output neuron (that predicts the spring height at time $t$).  Hyperbolic tangent activation was used at each hidden layer.  For the training problems, we used two types of terms: ODE-residual and data-fitting terms.  The times at which the ODE-residual terms were defined were 30 evenly spaced points over $[0,1]$.  The times at which the data-fitting terms were defined were 10 evenly spaced points over $[0,0.4]$.  The runs for \texttt{Adam(unc)} only considered an objective function where the terms in (\ref{eq.const}) were combined with a weight of $10^{-4}$ on the average ODE-residual.  The runs for the remaining solvers considered the same objective and included hard constraints for the ODE residual at times $\{\tfrac{4}{29},\tfrac{12}{29},\tfrac{21}{29}\}$, i.e., 3 constraints.  For all algorithms, we ran a ``full-batch'' version (i.e., with exact objective gradients employed) and a ``mini-batch'' version, where in each iteration of the latter version only half of the ODE-residual data points were used.  We employed the same two fixed learning rates (i.e., value of $\alpha$ in Algorithms~\ref{alg.heavy-ball} and~\ref{alg.padam_sqp}) for each algorithm: $5 \times 10^{-4}$ and $1 \times 10^{-4}$.  For the other step-size-related parameters we chose $\rho_k = 1$ and $h_k = 1$ for all $k \in \N{}$.

\begin{figure}[ht]
  \begin{center}
  \includegraphics[width=0.245\textwidth]{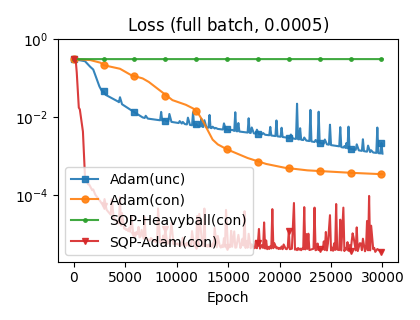}%
  \includegraphics[width=0.245\textwidth]{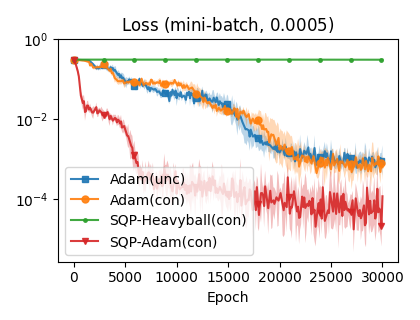}%
  \includegraphics[width=0.245\textwidth]{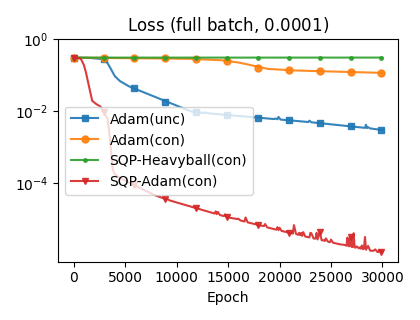}%
  \includegraphics[width=0.245\textwidth]{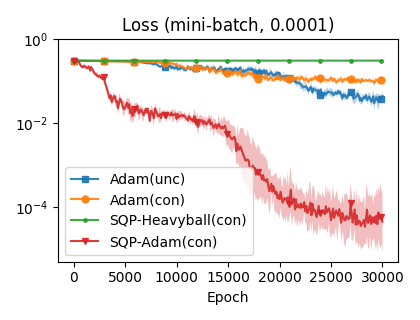}
  \caption{1D Spring losses over epochs. For the mini-batch runs, the solid lines indicate means over 5 runs while the shaded regions indicate values within one standard deviation of the means.}
  \label{fig:spring objective}
\end{center}
\end{figure}

\begin{figure}[ht]
  \begin{center}
  \includegraphics[width=0.245\textwidth]{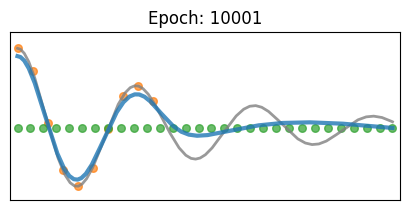}%
  \includegraphics[width=0.245\textwidth]{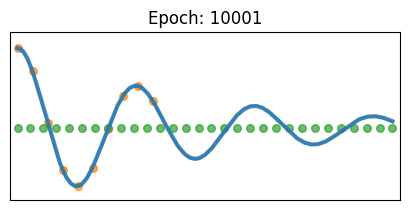}%
  \includegraphics[width=0.245\textwidth]{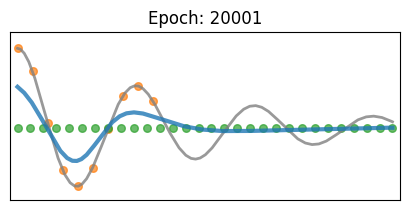}%
  \includegraphics[width=0.245\textwidth]{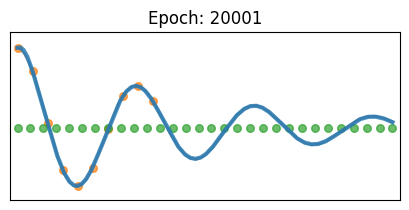}
  \caption{
    Predicted trajectories. Left to right: \texttt{Adam(unc)} (mini-batch, $\alpha_k=0.0005$), \texttt{SQP-Adam(con)} (mini-batch, 0.0005), \texttt{Adam(unc)} (mini-batch, 0.0001), and \texttt{SQP-Adam(con)} (mini-batch, 0.0001). Axes are time $t \in [0,1]$ (horizontal) and true/predicted $u(t)$ (vertical). Green dots indicate times at which the ODE-residual terms were defined for the objective function; orange dots indicate data-fitting values; the gray line indicates the true solution; and the blue line indicates the predicted solution.  Code from \citep{harmonic.oscillator} (available under the MIT License) was used to generate the plots.
  }
  \label{fig:spring solution}
\end{center}
\end{figure}

Results are provided in Figures~\ref{fig:spring objective} and \ref{fig:spring solution}.  The plots in Figure~\ref{fig:spring objective} show that \texttt{SQP-Adam(con)} yielded lower objective values (loss) more quickly and achieved better accuracy (i.e., lower mean-squared error) after the training budget expired.  They also show that \texttt{SQP-Adam(con)} achieved more comparable results for the two learning rates, whereas the other algorithms performed worse for the smaller learning rate. Our results here demonstrate that \texttt{SQP-Adam(con)} requires less hyperparameter tuning. The plots in Figure~\ref{fig:spring solution} show that the difference in performance can be seen clearly in the predictions that one obtains.

\subsection{Chemical engineering problem}\label{sec.chemistry}
This problem  models the reaction system of 1-butene isomerization when cracked on an acidic zeolite~\citep{Gupta2016}. The system is reformulated as an ordinary linear differential equation by scaling the kinetic parameters of the true model. The ODE is
\bequationNN
  \frac{du(t)}{dt} = \bbmatrix 
    -(c^{(1)} + c^{(2)} + c^{(4)}) u^{(1)}(t) + c^{(3)} u^{(3)}(t) + c^{(5)} u^{(4)}(t) \\
    2c^{(1)}u^{(1)}(t)\\
    c^{(2)}u^{(1)}(t) - c^{(3)} u^{(3)}(t)\\
    c^{(4)}u^{(1)}(t) - c^{(5)} u^{(4)}(t) \ebmatrix,
\eequationNN
where $c = [4.283, 1.191, 5.743, 10.219, 1.535]^T$.  Our aim was to train a neural network with the known ODE and mass-balance condition (namely, $\frac{du^{(1)}(t)}{dt} + 0.5 \frac{du^{(2)}(t)}{dt} + \frac{du^{(3)}(t)}{dt} + \frac{du^{(4)}(t)}{dt} = 0$) over various initial conditions near a nominal initial condition, where the nominal one is $u_0 = [14.5467,\  16.335,\  25.947,\  23.525]$.  In this manner, the trained network can be used to predict $u(t)$ at any $t$ (for which we use the range $t \in [0,10]$) for any initial condition near the nominal one.

\begin{figure}[ht]
  \begin{center}
  \includegraphics[width=0.245\textwidth]{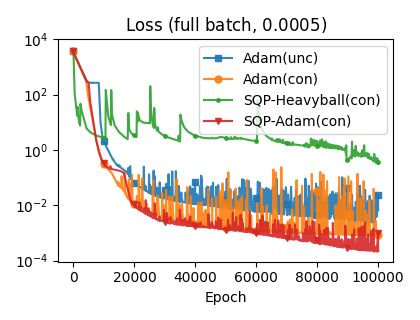}%
  \includegraphics[width=0.245\textwidth]{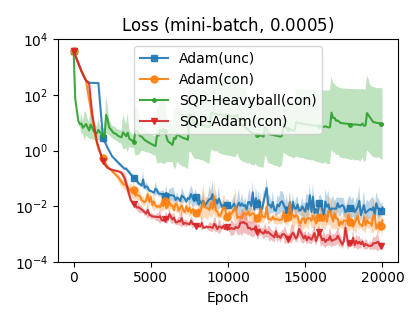}%
  \includegraphics[width=0.245\textwidth]{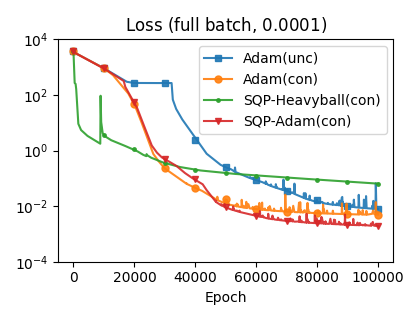}%
  \includegraphics[width=0.245\textwidth]{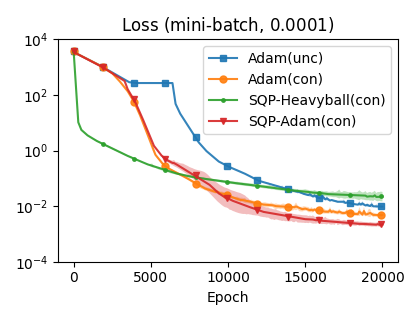}
  \caption{Chemical engineering problem losses over epochs. For the mini-batch runs, solid lines indicate means over 5 runs while the shaded regions indicate values within one standard deviation of the means.}
  \label{fig:chemical objective}
\end{center}
\end{figure}

We used a fully connected neural network with 5 input neurons (corresponding to the initial condition in $\R{4}$ and $t \in \R{}$), 3 hidden layers with 64 neurons each and hyperbolic tangent activation, and 4 output neurons (corresponding to $u(t) \in \R{4}$).  The training problems involved three objective terms: ODE-residual (weighted by $10^{-2}$), mass-balance (weighted by $10^{-2}$), and data-fitting (weighted by~$1$) terms.  Training data was generated by solving the ODE over 1000 initial conditions (of the form $u_0 + \xi$, where $\xi$ was a random vector with each element drawn from a uniform distribution over $[-1,1]$) using \texttt{odeint} from the \texttt{scipy} library \citep{2020SciPy-NMeth} (BSD licensed).  Specifically, solution values were obtained over 64 evenly spaced times in $[0,10]$, which over the 1000 initial conditions led to 64000 training points.  The ODE-residual and mass-balance terms involved all 64000 training points, whereas the data-fitting term involved only 20\% of these points chosen at random with equal probability.  The runs for \texttt{Adam(unc)} used only these objective terms, whereas the runs for \texttt{SQP-Heavyball(con)}, \texttt{SQP-Adam(con)}, and \texttt{Adam(con)} also considered 10 constraints on mass-balance residuals, the points for which were chosen uniformly at random over all initial conditions and times.  The mini-batch size was 20\% of all samples.  We tested learning rates for all algorithms: $5 \times 10^{-4}$ and $1 \times 10^{-4}$. For the other step-size-related
parameters we chose $\rho_k=0.5$ and $h_k=1$ for all $k \in \mathbb{N}$. The results in Figures~\ref{fig:chemical objective} 
show that \texttt{SQP-Adam(con)} performed best.

\subsection{1D Burgers' equation}

Burgers' equation is a PDE often used to describe the behavior of certain types of nonlinear waves \citep{Wang2021,Negiar2023}.  With respect to a spatial domain $[0,1]$, time domain $[0,1]$, and viscosity parameter $\nu = 0.01$, we used the equation, initial condition, and (periodic) boundary condition
\begin{align*}
  \tfrac{\partial u(x,t)}{\partial t} + u(x,t) \tfrac{\partial u(x,t)}{\partial x} &= \nu\tfrac{\partial ^{2} u(x,t)}{\partial x^{2}}, && x \in (0,1),\ t \in [0,1]; \\
  u(x,0) &= u_{0}(x), && x \in [0,1]; \\
  u(x,t) &= u(x+1,t), && x \in [0,1],\ t \in [0,1].
\end{align*}
Our aim was to train a neural network with the known PDE and boundary condition over various initial conditions near a nominal initial condition.  In this manner, for any $(x,t)$ and initial condition near the nominal one, the trained network can predict $u(x,t)$.

\begin{figure}[ht]
  \centering
  \includegraphics[width=0.245\textwidth]{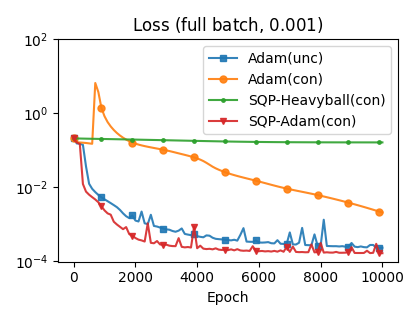}%
  \includegraphics[width=0.245\textwidth]{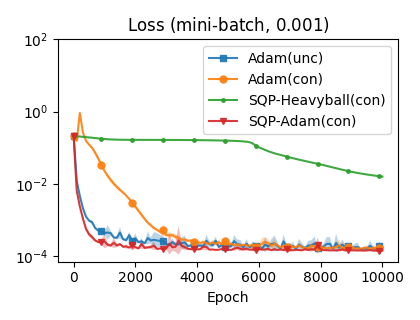}%
  \includegraphics[width=0.245\textwidth]{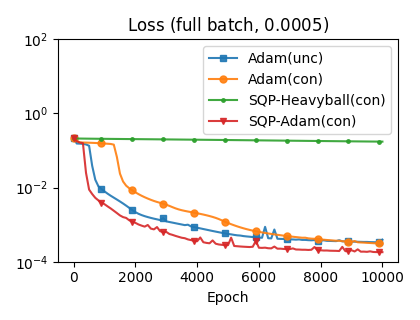}%
  \includegraphics[width=0.245\textwidth]{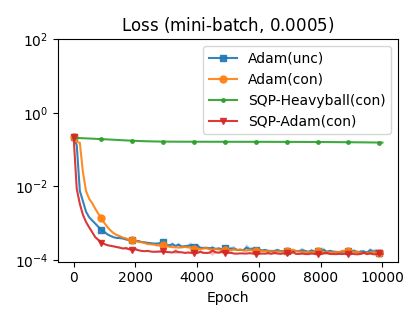}
  \caption{Burgers' losses over epochs. For mini-batch runs, solid lines indicate means over 5 runs while the shaded regions (not very visible) indicate values within one standard deviation of the means.}
  \label{fig:burgers objective}
\end{figure}

We used a fully-connected neural network with 34 input neurons (corresponding to $x$, $t$, and a discretization of $u_0$ over 32 evenly spaced points), 3 hidden layers with 64 neurons each and hyperbolic tangent activation, and 1 output neuron (corresponding to $u(x,t)$).  The training problems involve three objective terms: PDE-residual (weighted by $10^{-3}$), boundary-residual (weighted by $10^{-3}$), and data-fitting (weighted by 1) terms.  Training data was generated by solving the PDE over 100 initial conditions (of the form $u_0(x) = \sin(2\pi x + \xi{\pi})$, where for each instance $\xi$ was chosen from a uniform distribution over $[0,0.2]$) using the \texttt{odeint} solver, as in the previous section.  Specifically, solution values were obtained over 32 evenly spaced points each in the spatial and time domains, which over the 100 initial conditions led to 102,400 training points.  For each initial condition, the PDE-residual and boundary-residual terms involved all relevant generated training points, whereas the data-fitting term involved only 200 points chosen at random with equal probability.  \texttt{Adam(unc)} used only these objective terms, whereas \texttt{SQP-Heavyball(con)}, \texttt{SQP-Adam(con)}, and \texttt{Adam(con)} also considered 10 constraints on PDE residuals, the points for which were chosen uniformly at random over all initial conditions and spatio-temporal points.  The mini-batch was 20\% of all samples.  We tested learning rates: $10^{-3}$ and $5\times10^{-4}$. For the other step-size-related parameters we chose $\rho_k = 1$ and $h_k = 1$ for all $k \in \N{}$.  One finds in Figure~\ref{fig:burgers objective} that the results obtained by \texttt{Adam(unc)} and \texttt{SQP-Adam(con)} were in fact comparable.  The performance by the projection-less Adam approach (\texttt{Adam(con)}) was inferior to these, and the performance by \texttt{SQP-Heavyball(con)} was poorer still.  Figure~\ref{fig:burgers prediction} shows that a prediction by the model obtained by \texttt{SQP-Adam(con)} is indeed close to the true solution.

\begin{minipage}{\linewidth}
  \centering
  \begin{minipage}[t]{0.42\linewidth}
      \begin{figure}[H]
        \centering
        \includegraphics[width=3cm]{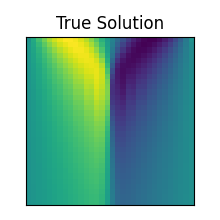}
        \includegraphics[width=3cm]{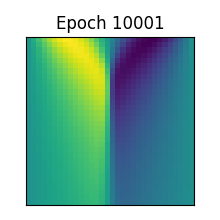}
        \caption{Burgers' true/predicted solutions for initial condition not seen in training. Predicted solution by \texttt{SQP-Adam(con)} (mini-batch, $\alpha_k=0.0005$).}
        \label{fig:burgers prediction}
      \end{figure}
  \end{minipage}
  \hspace{0.03\linewidth}
  \begin{minipage}[t]{0.51\linewidth}
      \begin{figure}[H]
        \centering
        \includegraphics[width=4.52cm]{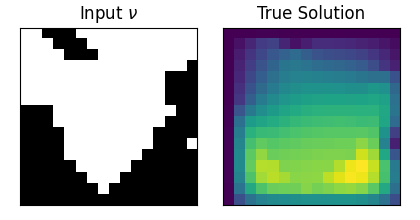}
        \includegraphics[width=2.48cm]{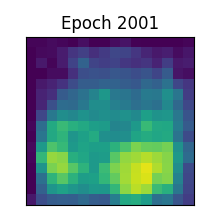}
        \caption{Darcy flow diffusion coefficient $\nu$ and true/predicted solution, where diffusion coefficient $\nu$ not seen in training. Predicted solution by \texttt{SQP-Adam(con)} (mini-batch, $\alpha_k=0.001$).}
        \label{fig:darcy prediction}
      \end{figure}
  \end{minipage}
\end{minipage}

\subsection{2D Darcy flow}

The steady-state 2D Darcy flow equations model the flow of a fluid through a porous medium \citep{Negiar2023,Takamoto2023}.  With respect to the spatial domain $[0,1]^2$, a forcing function $f$ (we use $f(x) = 1$ for all $x \in (0,1)^2$ in our experiments), and a diffusion coefficient $\nu$, we used
\begin{align*}
  -\nabla \cdot (\nu(x) \nabla u(x)) &= f(x), && x \in (0,1)^2; \\ u(x) &= 0, && x \in \partial [0,1]^2.
\end{align*}
Our aim was to train a neural network with the known PDE and boundary condition over various diffusion coefficients such that, for any $x \in [0,1]^2$ and diffusion coefficient $\nu$, it could be used to predict $u(x)$.

\begin{figure}[ht]
    \centering
    \includegraphics[width=0.245\textwidth]{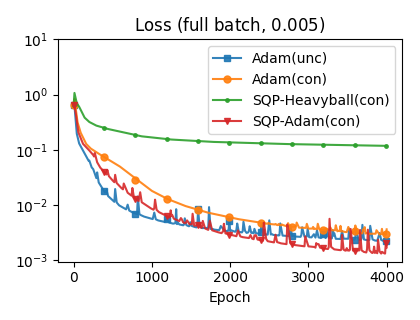}%
    \includegraphics[width=0.245\textwidth]{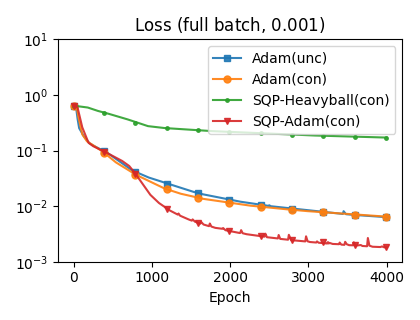}%
    \includegraphics[width=0.245\textwidth]{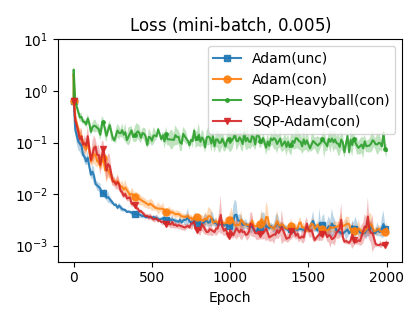}%
    \includegraphics[width=0.245\textwidth]{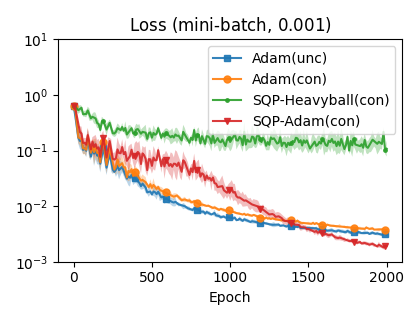}
  \caption{Darcy flow losses over epochs. For mini-batch runs, solid lines indicate means over 5 runs while the shaded regions indicate values within one standard deviation of the means.}
  \label{fig:darcy objective}
\end{figure}

We used the Fourier Neural Operator (FNO) architecture imported from the \emph{neuralop} library \citep{kovachki2021neural, li2021fourier} (MIT License). The inputs were given in three channels, one for $\nu$, one for a horizontal position embedding, and one for a vertical position embedding.  Each channel had dimension $16\times 16$.  We used 4 hidden layers (the default). The output was a single channel of dimension $16\times 16$ (corresponding to $u(x)$).  The training problems involve three objective terms: PDE-residual (weighted by $10^{-3}$), boundary-residual (weighted by $10^{-3}$), and data-fitting (weighted by 1) terms. We use the \emph{neuralop} package \cite{kossaifi2025librarylearningneuraloperators} to generate 100 $\nu$ values and their corresponding solutions for training. Specifically, we first generate 1000 samples using the default settings of \emph{neuralop} and then select 100 with similar $\nu$ values. For each $\nu$ value, the PDE-residual and boundary-residual terms involved all relevant generated training points, whereas the data-fitting term involved only 20\% of the points chosen at random with equal probability.  The runs for \texttt{Adam(unc)} used only these objective terms, whereas the runs for \texttt{SQP-Heavyball(con)}, \texttt{SQP-Adam(con)}, and \texttt{Adam(con)} considered the same objective in addition to 50 constraints on PDE residuals, the points for which were chosen uniformly at random over all initial conditions and spatial points.  We ran full-batch and mini-batch settings, where the mini-batch was dictated by 20\% of the $\nu$ values.  We tested using the same learning rates for all algorithms: $5\times10^{-3}$ and $1\times10^{-3}$.  For the other step-size-related parameters we chose $\rho_k = 1$ (\texttt{SQP-Heavyball(con)}), $\rho_k = 0.5$ (\texttt{SQP-Adam(con)}), and $h_k = 1$ for all $k \in \N{}$.  The results in~\ref{fig:darcy objective} show strong relative performance by \texttt{SQP-Adam(con)}.  Figure~\ref{fig:darcy prediction} shows that a prediction by the model obtained by \texttt{SQP-Adam(con)} is close to the true solution.

\section{Conclusion}\label{sec.conclusion}

We proposed two stochastic diagonal-scaling methods for nonlinear equality constrained optimization, and provided convergence guarantees for each approach.  We also demonstrated the algorithms in the context of informed supervised learning.  The methods' per-iteration costs are comparable to an unconstrained (soft-constrained) approach that also uses diagonal scaling.  The numerical experiments reveal practical benefits of the proposed schemes, which we conjecture would also be witnessed when training larger and more sophisticated neural networks for informed learning.

\bibliographystyle{plain}
\bibliography{references}

\end{document}